\documentclass[11pt]{article}

\usepackage{latexsym,amsmath,amsfonts,amssymb,amstext,theorem,euscript,array,enumerate,mathrsfs,color}

\newtheorem{Theorem}{Theorem}[section]
\newtheorem{Definition}{Definition}[section]
\newtheorem{Proposition}{Proposition}[section]

\newtheorem{Lemma}{Lemma}[section]

\newtheorem{Example}{Example}[section]
\numberwithin{equation}{section}

\def \r{\mathfrak{r}}
\def \hX{\hat{X}}
\def \hY{\hat{Y}}
\def \hZ{\hat{Z}}
\def \hU{\hat{U}}

\def \R{\mathbb{R}}

\def \E{\mathbb{E}}

\def \P{\mathbb{P}}

\def \l{{\mathfrak{g}}}

\def \F{{\cal F}}

\def \Pc{{\cal P}}

\newcommand{\tr}{\text{tr}}
\def\e{\varepsilon}

\def \ep{\hbox{ }\hfill$\Box$}
\def \pr{\textbf{Proof.}}
\allowdisplaybreaks

\begin{document}

\title{Ergodic BSDEs with Multiplicative and Degenerate Noise}

\author{ 
Giuseppina GUATTERI\footnote{Department of Mathematics, Politecnico di Milano, Piazza Leonardo da Vinci 32, 20133 Milano, Italy, {\tt giuseppina.guatteri@polimi.it}.}
\qquad\quad
Gianmario TESSITORE\footnote{Department of Mathematics and Applications, University of Milano-Bicocca, Via R. Cozzi 53 - Building U5, 20125 Milano, Italy, {\tt gianmario.tessitore@unimib.it}.}} 

\maketitle

\begin{abstract}
In this paper we study an Ergodic Markovian BSDE  involving a forward process $X$ that solves an infinite dimensional forward stochastic evolution equation with multiplicative  and possibly degenerate diffusion coefficient.
A concavity assumption on the driver allows us to   avoid the typical quantitative conditions relating the dissipativity of the forward equation and the Lipschitz constant of the driver.
Although the degeneracy of the noise has to be of a suitable type we can give a stochastic representation of a large class of Ergodic HJB equations; morever  our general results can be applied to get the synthesis of the optimal feedback law  in relevant 
examples of ergodic  control problems for SPDEs.

\end{abstract}

\noindent {\bf Keywords:} Ergodic control; infinite dimensional SDEs; BSDEs; Multiplicative Noise

\noindent {\bf 2010 Mathematics Subject Classification:} 60H15, 60H30, 37A50.

\date{}

\section{Introduction}
In this paper we study the following BSDE of ergodic type
$$
Y^x_t= Y^x_T + \int_t^T [\widehat{\psi}(X^x_s,  Z^x_s,U^x_s )-\lambda] \, ds - \int_t^T  Z^x_s \, dW^1_s-\int_t^T  U^x_s \, dW^2_s, \qquad   0 \leq t \leq T < \infty, 
$$
where the processes $(Y^x,Z^x,U^x)$ and the constant $\lambda$ are the unknowns of the above equation while the diffusion $X$ is the (mild) solution of the infinite dimensional (forward) SDE:
 \begin{equation*}\left \{ \begin{array}{l}
 d X^{x}_s \ = \ A X^{x}_sds + F( X^{x}_s)ds + Q G(X^{x}_s)d W^1_s + Dd W^2_s , \\  X^{t,x}_t \ = \ x.
\end{array}\right.
\end{equation*}
In the above equation $X$ takes values in an Hilbert space $H$ and $W^1$, $W^2$ are independent cylindrical Wiener processes (see \textbf{(A.1)}-\textbf{(A.6)} in Section \ref{forward} and \textbf{(B.1)} in Section \ref{backward} for precise description of the other terms). We just stress that we will assume that $G(x)$ is invertible for all $x\in H$ while $Q$ and $D$ will be general, possibly degenerate, linear operators.

Ergodic BSDEs have been introduced in \cite{FuhHuTess} in relation to optimal stochastic ergodic control problems  and as a tool to study the asymptotic behaviour of parabolic HJB equations and consequently to give a stochastic representation to the limit semilinear elliptic PDEs (see equation \eqref{HJB} below).

 In \cite{FuhHuTess} the same class of BSDEs have been introduced, already  in an infinite dimensional framework, but only in the case in which the noise coefficient was  constant ($Q=0$ in our notation). Successive works, see \cite{HuMadRic} and \cite{DebHuTess} weakened the assumptions and refined the results in the same \textit{additive noise} case. Then in \cite{Richou20092945}, in a finite dimensional framework, the case of `multiplicative noise ($Q \neq 0$ and $G$ depending on $x$ in our notation) is treated under quantitative conditions relating the dissipativity constant of the forward equation to the Lipscitz norm of $\widehat{\psi}$ with respect to $Z$. Afterwards, in \cite{Madec20151821}, still in finite dimensions, such quantitative assumptions are dropped in the case of a non degenerate and bounded diffusion coefficient ($Q=I$ and $G$ bounded and invertible in our notation) by a careful use of smoothing properties of the Kolmogorov semigroup associated to the non-degenerate underlying diffusion $X$. Finally in  
 \cite{HuLemmonier2018} the result is extended to the case of non degenerate but unbounded (linearly growing)  diffusion coefficients  ($Q=I$ and $G$ invertible and linearly growing in our notation). 
To complete the picture we mention, \cite{Cohen20134138}, \cite{CoFuhPham},  \cite{CossoGuatteriTessitore} and \cite{Hu2018}
where Ergodic BSDEs are studied in various frameworks different from the present one: namely, respectively  when they are driven by a Markov chain, in the context (see \cite{KharroubiPham}) of  randomized control problems and BSDEs with constraints on the martingale term both in finite and in infinite dimensions and finally in the context of $G$- expectations theory.

In this paper we propose an alternative approach that works well in the infinite dimensional case and allows to consider degenerate multiplicative noise ($Q$ in general non invertible and $G$ bounded invertible but depending on $x$). On the other side we have to assume that $\widehat{\psi}$  has the form:
{    $$\widehat{{\psi}}(x,z ,u):=\psi(x,zG^{-1}(x),u)$$ where $\psi$  is Lipschitz and \textit{concave} function with respect to $(z,u)$.
}
 Although not standard, our assumptions allow to give a stochastic representation of a relevant class of Ergodic HJB equations in Hilbert spaces (see Section \ref{sec_HJB}) and of ergodic stochastic control problems for SPDEs (see Example \ref{esempio_bordo} and Example \ref{esempioaltro}). Notice that $\psi$ defined above is exactly the function that naturally appears in the related HJB equation and in the applications to ergodic control.

As in all the literature devoted to the problem the main point is to prove a uniform gradient estimate (independent on $\alpha$) for  $v^{\alpha}(x):=Y^{\alpha,x}$ where  $(Y^{\alpha,x},Z^{\alpha,x},U^{\alpha,x})$  is the solution of the discounted BSDE with infinite horizon:
$$Y^{\alpha,x}_t= Y^{\alpha,x}_T + \int_t^T [\widehat{\psi}(X^x_s,  Z^{\alpha,x}_s,U^{\alpha,x}_s )-\alpha Y^{\alpha,x}_s] \, ds - \int_t^T  Z^{\alpha,x}_s \, dW^1_s-\int_t^T  U^{\alpha,x}_s \, dW^2_s, \; \quad   0 \!\leq\! t \!\leq \! T \!<\! \infty, 
$$
Such estimate can be obtained by a change of probability argument when the noise is additive (see \cite{FuhHuTess}), by energy type estimates under quantitative assumptions on the exponential decay of the forward equation (see \cite{Richou20092945}) or by regularizing properties of the Kolmogorov semigroup when the noise in multiplicative but non degenerate (see \cite{HuLemmonier2018} and \cite{Madec20151821}).

Here we exploit concavity of $\psi$ to introduce an auxiliary control problem and eventually obtain the gradient estimate using a decay estimate on the difference between states starting from different initial conditions,  see Assumption \textbf{(A.6)} and, in particular, requirement \eqref{stimadiffSol}. We stress the fact that the estimate in \eqref{stimadiffSol} is only in mean and not uniform (with respect to the stochastic parameter) as in the additive noise case.
Moreover, as we show in Proposition \ref{propdiss}, Assumption \textbf{(A.6)} is verified if we impose a \textit{joint dissipativity} condition on the  coefficients, see Assumption \textbf{(A.7)}. As a matter of fact, in this case, the stronger formulation in which $L^2$  replaces $L^1$ norm holds. On the other side \textbf{(A.6)}  allows to cover a wider class of interesting examples, see for instance Example \ref{esempio_bordo} in which  Assumption \textbf{(A.7)} does not seem to hold. 

The structure of the paper in the following: in Section
\ref{notation} we introduce the function spaces that will be used in the following,  Section   \ref{forward} is devoted to the infinite dimensional forward equation; in particular we state and discuss the  key stability assumption \textbf{(A.6)}. In Section \ref{backward} we present the main contribution of this work introducing the auxiliary control problem, proving the gradient estimate and the consequent existence of the solution to the ergodic BSDEs.  In Section \ref{sec_HJB} we relate our ergodic BSDE to a semilinear PDE in infinite dimensional spaces (the ergodic HJB equation). In Section \ref{sub-diff} we discute the regularity of the solution of the ergodic BSDE, in particular we state that under quantitative conditions on the dissipativity of the forward equation similar to the ones assumed in \cite{Richou20092945}, when all coefficients are differentiable then the solution of the ergodic BSDE is differentiable with respect to the initial data as well. The proof of such result adapts a similar argument in  \cite{HuTess} and is rather technical,  we have postponed it in the Appendix In Section \ref{sec_control} we use our ergodic BSDE to obtain an optimal ergodic control problem (that is with cost depending only on the asymptotic behaviour of the state) for an infinite dimensional equation. We close, see Section \ref{sub-ex}, by two examples of controlled SPDEs to which our results can be applied. In both we consider a stochastic  heat equation in one dimension with additive white noise. In the first, Example \ref{esempio_bordo} the system is controlled through one Dirichlet boundary condition (on which multiplicative noise also acts) while, in the second one, Example \ref{esempioaltro}, the control enters the system through a finite dimensional process that affects  the coefficients of the SPDE. In this last case we also give conditions guaranteeing differentiability of the  related solution to the Ergodic BSDE.

\section{General notation}\label{notation}

Let $\Xi$, $H$ and $U$ be real separable Hilbert spaces. In the sequel, we use the notations $|\cdot|_\Xi$, $|\cdot|_H$ and $|\cdot |_U$ to denote the norms on $\Xi$, $H$ and $U$ respectively; if no confusion arises, we simply write $|\cdot|$. We use similar notation for the scalar products. We denote the dual spaces of $\Xi$, $H$ and $U$ by $\Xi^*$, $H^*$,  and $U^*$ respectively. We also denote by $L(H,H)$  the space of bounded linear operators from $H$ to $H$, endowed with the operator norm. Moreover, we denote by $L_2(\Xi,H)$ the space of Hilbert-Schmidt operators from $\Xi$ to $H$. Finally,  a map $f: H \to \Xi$ is said to belong to the class $\mathcal{G}^1(H,\Xi)$ if it is continuous and Gateaux differentiable with  directional derivative $\nabla _x  f(x)h$ in $(x,h) \in H \times H $   and we denote by $\mathcal{B}(\Lambda)$ the Borel
$\sigma$-algebra of any topological space $\Lambda$.

Given a complete probability space $(\Omega, \mathcal{F}, \mathbb{P})$ together with a filtration $(\mathcal{F}_t)_{t\geq0}$ (satisfying the usual conditions of $\P$-completeness and right-continuity) and  an arbitrary real separable Hilbert space $V$  we define the following classes of processes for fixed $0\leq t\leq T $ and  $p\geq 1$:
\begin{itemize}
\item   $L^p_\Pc (\Omega\times [t,T];V)$ denotes the set of (equivalence classes) of $(\mathcal{F}_s)$-predictable processes $Y \in L^p (\Omega\times
[t,T];V)$ such that the following norm is finite:
\[
|Y|_p \ = \ \bigg(\E \int_t^T |Y_s|^p \, ds\bigg)^{1/p}
\]
\item $L^{p, loc}_{\mathcal{P}}(\Omega\times [0,+\infty[;V)$ denotes the set of processes defined on $\mathbb{R}^+$, whose restriction to an arbitrary time interval $[0,T]$ belongs to $L^p_\Pc (\Omega\times [0,T];V) $.  
   \item $L^p_\Pc(\Omega;C([t,T];V))$
     denotes the set of $(\mathcal{F}_s)$-predictable processes $Y$ on $[t,T]$ with continuous paths in $V$, such
    that the norm
    \[  \|Y\|_p \ = \ \big(\E \sup _{s \in [t,T]} |Y_s|^p\big)^{1/p}\]
    is finite. The elements of $L^p_{\mathcal{P}}(\Omega;C([t,T];V))$
    are identified up to indistinguishability.
    \item $L^{p, loc}_\Pc (\Omega;{C} ([0,+\infty[;V))$ denotes the set of processes defined on $\mathbb{R}^+$, whose restriction to an arbitrary time interval $[0,T]$ belongs to $L^p_\Pc(\Omega;C([0,T];V))$. 
\end{itemize}
We consider on the 
probability space $(\Omega,\F,\P)$ two independent cylindrical Wiener processes $ W^1 = (W^1_t)_{t\geq 0}$  with
values in $\Xi$ and  $ W^2 = (W^2_t)_{t\geq 0}$  with
values in $H$.
By $(\F_t)_{t\geq0}$, we denote the natural filtration
of $( W^1, W^2) $, augmented with the family ${\mathcal N}$ of
$\P$-null sets of $\F$. The filtration
$(\F_t)$ satisfies the usual conditions of right-continuity and $\P$-completeness.

\section{Forward equation}\label{forward}

Given $x\in H$ and a uniformly bounded process $\mathfrak{g}$ with values in $H$, we consider the  stochastic differential equation for $t \geq 0$
\begin{equation}\label{State1}
d X^{x,\l}_t \ = \ A X^{x,\l}_tdt + F( X^{x,\l}_t)dt + Q G(X^{x,\l}_t)d W^1_t  + D d W^2_t+  \l (t)\, dt, \qquad  X^{x,\l}_0 \ = \ x.
\end{equation}
On the coefficients $A$, $F$, $G$, $Q$, $D$   we impose the following assumptions.

\begin{itemize}
\item[({\bf A.1})] $A\colon\mathcal{D}(A)\subset H\to H$ is a linear, possibly unbounded operator   generating a $C_0$ semigroup  $\{e^{tA}\}_{t\geq 0}$.

\item[({\bf A.2})] $F\colon H\to H$ is continuous and there exists  
%
 $L_F>0$ such that
\[
|F(x) - F(x')|_H \ \leq \ L_F  |x - x'|_H,
\]
for all $x,x'\in H$.

\item[({\bf A.3})] $G\colon H \rightarrow  L(\Xi)$  is a bounded Lipschitz map. Moreover, for every $x \in H$, $G(x)$ is invertible. Thus
 there exists  three    positive constants $L_G$, $M_{G}$ and $M_{G^{-1}}$
 such that for all  $x,x'\in H$:
\[
|G(x)|_{L(\Xi)}\ \leq {M_G } \qquad 
|G(x) - G(x')|_{L(\Xi)}\ \leq {L_G }  |x - x'|_H,\qquad  \big|G^{-1} (x)\big|_{L(\Xi)} \leq {M_{G^{-1}} }
\]
We notice that the above yields Lipschitzianity of $G^{-1}$, namely :
 \[
|G^{-1}(x) - G^{-1}(x')]|_{L(\Xi)}\ \leq { M^2_{G^{-1}}L_{G}} \,  |x - x'|_H,
\]


\item[({\bf A.4})]  $Q$ is an Hilbert-Schmidt operator from $\Xi$ to $H$.
\item[({\bf A.5})]  $D$ is a linear and bounded operator from $H$ to $H$ and there exist  constants $L>0$ and $\gamma \in [0,\frac{1}{2}[$:
\begin{align}
|e^{sA}D|_{L_2(H)} \leq L
\left(s^{-\gamma}\wedge 1\right),\quad \forall s\geq 0.\end{align}

\end{itemize}
 \begin{Proposition}\label{Prop-esunstate} 
Under $({\bf  A.1} -- {\bf  A.5} )$,  for any $x\in H$ and any $\l$ bounded and progressive measurable process with values in $H$, there exists a unique (up to indistinguishability) process $X^{x,\l}=( X_t^{x,\l})_{t\geq0}$ that belongs to $L^{p,loc}_\Pc (\Omega;C([0,+\infty[;H))$ for all $p \geq 1$ and is a mild solution of \eqref{State1}, that is it satisfies for every  $ t\!\geq\!0$, $\P \text{-a.s.}$:
\begin{align*}\label{mild}
 X_t^{x,\l}   &=  e^{tA} x  + \int_0^t e^{(t-s)A} F( X_s^{x,\l}) \, ds + \int_0^t e^{(t-s)A} \l(s) \, ds +  \int_0^t e^{(t-s)A}QG( X_s^{x,\l}) \, d W^1_s \\ & + \int_0^t e^{(t-s)A}D \, d W^2_ s. \ 
\end{align*}
Moreover there exists a positive constant $\kappa_{\l, T}$  such that
\begin{equation}\label{stimaSolgamma}
{\E}|{ X}^{x,\l}_t| ^2\ \leq \ \kappa_{\l ,T}(1+ |x|^2), \qquad \qquad \forall t \in [0,T] \text{ and  } x \in H.
\end{equation}

\end{Proposition}
Our main result will be obtained  under the following exponential stability in $L^1$ norm requirement. We stress the fact that such assumption is much weaker in comparison with the uniform decay holding when noise is addittive (see \cite{FuhHuTess}).

\begin{itemize}
\item[({\bf A.6})]
 There exist positive constants $\kappa_{\l}$,  $\kappa$ and $\mu$ such that
\begin{equation}\label{stimaSol}
\sup_{t\geq 0}{\E}|{ X}^{x,\l}_t|  \leq \kappa_{\l}(1+ |x|);
\end{equation}
 \begin{equation}\label{stimadiffSol}
{\E}|{ X}^{x,\l}_t - { X}^{x',\l}_t| \ \leq \kappa   e^{-\mu t}|x- x'|;
\end{equation}
 for any $x,x' \in H$ and  for all $ t  \geq 0$.
\end{itemize}
Below we show that hypothesis $({\bf A.6})$ (as a matter of fact the stronger condition obtained replacing $L^1$ norm by $L^2$ norm) is verified under the usual {\em joint dissipative condition}  $({\bf A.7})$ (see \cite{DPZab2}). We { have preferred} to keep the weaker, but less intrinsic,  form  $({\bf A.6})$ since it allows to cover a wider class of examples, see for instance Example \ref{esempio_bordo}
\begin{itemize}
\item[({\bf A.7})]{ - \em Joint dissipative conditions}

\noindent $A$ is dissipative i.e. $<Ax,x>\,\leq  \rho |x|^2 $, for all $x\in \mathcal{D}(A),$  and for some $\rho\in \mathbb{R}$, moreover
there exists $\mu>0$ such that for all $x, x'\in D(A)$:
\begin{equation}\label{condiss}
2\langle  A (x-x')+ F(x) -  F(x'), x - x' \rangle_H  + || Q [G(x) - G(x')] ||^2_{L_2(\Xi,H)}\ \leq \ -\mu|x - x'|_H^2,
\end{equation}
Notice that, by adding a suitable constant to $F$ and subtracting it from $A$ we can always assume that $\rho$ above is strictly negative.
\end{itemize}
Indeed we have that following holds
\begin{Proposition}\label{propdiss}
Assume $({\bf A.1}--{\bf A.5})$ and
$({\bf A.7})$  then the following estimates hold for the solution $ { X}^{x,\l}$ of  equation \eqref{State1}:
\begin{equation}\label{stimaSolqua}
\sup_{t\geq 0}{\E}|{ X}^{x,\l}_t|^2  \leq \kappa_{\l}(1+ |x|^2);
\end{equation}
 \begin{equation}\label{stimadiffSolqua}
{\E}|{ X}^{x,\l}_t - { X}^{x',\l}_t|^2 \ \leq   e^{-\mu t}|x- x'|^2;
\end{equation}
 for any $x,x' \in H$ and  for all $ t  \geq 0$.
In particular, hypothesis $({\bf A.6})$ is verified.
 \end{Proposition}
\textbf{Proof.} 

The proof of these estimates follows rather standard arguments, for the reader's convenience we give some details in particular on the way infinite dimensionality of the state space can be handled.

Let $ V_s =\displaystyle \int_0^s e^{(s-r)A} Dd W^2_r  + \int_0^s e^{(s-r)A}  \l(r) dr$ and $\chi^x_t: = X^x_t - V_s$, then
\begin{equation}\label{StateZ}
d \chi^x_t \ = \ A \chi^x_tdt + F( X^x_t ) dt+ Q G(X^x_t ) d W^1_t \qquad  \chi^x_0 \ = \ x.
\end{equation}

\noindent For  any $ n \in \mathbb{N}$ consider $ J(n,A) :=(nI -A)^{-1}$ and define $$X^{n,x}_s:=J(n,A)X^x_s,\quad \chi^{n,x}_t :=
J(n,A)\chi^{x}_t=J(n,A) X^x_t -  J(n,A) V_t $$


\noindent It is well known that $ \sup_{n \geq 0}|J(n,A)| _{L(H)}  < \infty $  and  $\lim_{n\to \infty} J(n,A) x = x, \ \forall x \in H$ with the obvious consequences on the $\mathbb{P}$-a.s and $L^p(\Omega)$ convergence of $X^{n,x}_s$ towards $X_s$ and $\chi^n_s$ towards $\chi_s$.

\noindent By easy computations $  \chi^n_t $ solves:
$$
d \chi^{n,x}_t = A \chi^{n,x}_tdt + F\left( X^{n,x}_t \right)dt + Q G\left( X^{n,x}_t\right)d W_t^1 \nonumber  + R^{n,x}_t \, dt + S^{n,x}_t \, d W_t^1,\qquad \chi^n_0  = \ J(n,A)x,$$
where 
\begin{equation*}
R^{n,x}_t= J(n,A)  F( X^x_t) -  F( X^{n,x}_t ), 
\qquad
S^{n,x}_t= J(n,A)  QG( X^x_t ) -  QG(  X^{n,x}_t ).
\end{equation*}
From  hypotheses $\bf{(A.2)}$ and  $\bf{(A.3)}$  we deduce that:
$$
|R^{n,x}_t|_H \leq C (1+ |X^x_t|_H), \qquad |S^{n,x}_t|_{L_2(\Xi,H)} \leq C (1+ |X^x_t|_H).$$
Moreover, we have that for all $t\geq 0$ and all $x \in H$:
\begin{equation}\label{LIMRn}
 \lim_{n\to + \infty} |R^{n,x}_t|^2 \to 0, \qquad \mathbb{P}-a.s.,
\end{equation}
and, by a dominated convergence argument on the computation of the Hilbert Schmidt norm, see also  \cite[Lemma 5.1]{GuaTess08}, we  have that for all $t\geq 0$ and all $x \in H$
\begin{equation}\label{LIMSn}
\lim_{n\to + \infty}  |S^{n,x}_t|^2_{L_2(\Xi,H)} \to 0, \qquad \mathbb{P}-a.s.
\end{equation}
We apply It\^o's formula to $ e^{\mu t}|\chi^{n,x}_t|^2$, and we add and subtract terms in order to apply the joint dissipativity condition in $\textbf{(A.7)}$
\begin{align*}
&e^{\mu t}|\chi^{n,x}_t|^2 = |x|^2+ 2 \int_0^t e^{\mu s} \langle \chi^{n,x}_s, \frac{\mu}{2} \chi^{n,x}_s + A_n \chi^{n,x}_s  \rangle _H \, ds \\ &
+2 \int_0^t e^{\mu s} \langle \chi^{n,x}_s, F(X^{n,x}_s) - F(J(n,A)V_s )\rangle _H \, ds \\& + 2 \int_0^t e^{\mu s} \langle \chi ^{n,x}_s, F(J(n,A)V_s ) \rangle _H \, ds   + 2 \int_0^t e^{\mu s} \langle \chi^{n,x}_s, Q G(X^{n,x}_s)d W^1_s  \rangle _H \\ &+
\int_0^t e^{\mu s}\ \tr\left[\left[G(X^{n,x}_s) - G\left( F(J(n,A)V_s \right) \right]^T Q^T Q\left[G(X^{n,x}_s) - G\left(F(J(n,A)V_s\right) \right]\right]\, ds  \\ &+ 2 \int_0^t e^{\mu s}\tr\left[G(X^{n,x}_s)^T  Q^T Q G\left(F(J(n,A)V_s\right) \right]\, ds \\& -  \int_0^t e^{\mu s}\tr \left[G\left( J(n,A)V_s \right)^T  Q^T Q G\left( J(n,A)V_s \right) \right]\, ds  \\  & +2 \int_0^t e^{\mu s} \langle \chi ^{n,x}_s, R^n_s \rangle _H \, ds +
 2 \int_0^t e^{\mu s} \langle \chi^{n,x}_s, S^n_sd W^1_s  \rangle _H   + \int_0^t e^{\mu s}\tr\left[S^n_s (S^n_s)^T \right]\, ds \\ & + \int_0^t e^{\mu s}\tr \left[G(X^{n,x}_s )^T  Q^T  S^n_s  + S^n_sQ G(X^{n,x}_s ) \right]\, ds\end{align*}
and by $\eqref{condiss}$:
\begin{align*}
&e^{\mu t}|\chi^{n,x}_t|^2 \leq  |x|^2+ 2  \int_0^t e^{\mu s} \langle \chi ^{n,x}_s, F(J(n,A)V_s)  \rangle _H \, ds  +
 2 \int_0^t e^{\mu s} \langle \chi^{n,x}_s, Q G(X^n_s )d W^1_s  \rangle _H \\& +
 2 \int_0^t e^{\mu s} \langle \chi^{n,x}_s, S^n_sd W^1_s  \rangle _H    +  \int_0^t e^{\mu s} |S^n_s|^2_{L_2(\Xi,H)}\, ds +  2 \int_0^t e^{\mu s} \langle \chi ^{n,x}_s, R^n_s \rangle _H \, ds  \\&  +\int_0^t e^{\mu s}\tr \left[G(X^{n,x}_s )^T  Q^T Q G\left( J(n,A)V_s\right) + G\left( J(n,A)V_s \right)^T  Q^T Q G(X^{n,x}_s ) \right ]\, ds \\&  + \int_0^t e^{\mu s}\tr [G(X^{n,x}_s )^T  Q^T  S^n_s  + S^n_sQ G(X^{n,x}_s)  ]\, ds.
\end{align*}
By $({\bf A.3})$ and $({\bf A.5})$ the definition of $S^n$  and  the estimate \eqref{stimaSolqua} we have that the stochastic integrals  are martingales, and 
$$
e^{\mu t}\E |\chi^{n,x}_t|^2 \leq  |x|^2+ \frac{\mu}{2}\E  \int_0^t e^{\mu s} |\chi ^{n,x}_s|^2_H \, ds + C \E  \int_0^t e^{\mu s}\left( |R^{n,x}_s|^2_H\, ds  + |S^{n,x}_s|^2_{L_2(\Xi,H)}+s^{ -2 \gamma}+1\right)\, ds, $$
where $C$ is a constant independent of $t$ and $n$.



 Limits \eqref{LIMRn}, \eqref{LIMSn} and the Dominated Convergence Theorem imply that for every $ t\geq 0$, $$\displaystyle\lim_{n\to \infty} \E  \int_0^t e^{\mu s}  |R^{n,x}_s|^2_H\, ds =0$$ and 
$$\displaystyle\lim_{n\to \infty} \E  \int_0^t e^{\mu s}  |S^{n,x}_s|^2_{L_2(\Xi,H)}\, ds =0,$$ 

Therefore, letting $n$ tend to $\infty$:
$$
\E |\chi^x_t|^2 \leq  |x|^2+ \frac{\mu}{2}  \int_0^t e^{-\mu (t-s)} \E |\chi^x_s|^2_H \, ds + C  \int_0^t e^{-\mu (t-s)}\left( s^{ -2 \gamma}+1\right)\, ds, $$
and 
$$
\sup_{s\leq t}\E|\chi^x_s|^2 \leq  |x|^2+ \frac{\mu}{2}   \sup_{s\leq t}\E|\chi^x_s|^2 \int_0^t e^{-\mu (t-s)}  \, ds + C_1 $$
where $C_1$ depends on $\mu$ and $\gamma$ but not on $t$. Thus we can conclude that 
$
\E|\chi^x_s|^2 \leq C_2(1+|x|^2)$, for all $s\geq 0$ and that, for all $t \geq 0$:
\begin{equation}
\E |X^x_t|^2  \leq 4 \left(\E |\chi^x_t|^2 +  \left| \int_0^t e^{(t-r)A} Dd W^2_r\right| ^2 +  \left| \int_0^t e^{(t-r)A}    \l(r)dr\right| ^2\right) \leq  C (1+ |x|^2).
\end{equation}
where the constant $C$ is independent from $t$ thanks to the dissipativity assumptions on $A$.

Estimate \eqref{stimadiffSolqua}  follows by the similar (and indeed easier arguments) applying It\^o formula to the difference $|X^{n,x}_t - X^{n,x'}_t|^2=|\chi^{n,x}_t-\chi^{n,x}_t|^2$ noticing that:
$$
d (\chi^x_t -\chi^{x'}_t)= A (\chi^x_t -\chi^{x'}_t)dt + [F( X^x_t )- F( X^{x'}_t) ]dt+ Q [G(X^x_t ) -G(X^{x'}_t )]d W^1_t \qquad  \chi^x_0-\chi^{x'}_0 =0.
$$
\ep

We end this section noticing that will be  mainly interested in the special case where $\l \equiv 0$:
\begin{equation}\label{State1null}
d X_t \ = \ A X_tdt + F( {X}_t)dt + Q G(X_t)d W^1_t  + D d W^2_t, \qquad  X^{x}_0 \ = \ x,
\end{equation}
and we will denote by $ X^{x}$ its solution through the whole paper.

\section{Ergodic BSDEs  }\label{backward}

In this section we study the following  equation:
\begin{equation}\label{ergodicbsde}
Y^x_t= Y^x_T + \int_t^T [\psi(X^x_s,  Z^x_sG^{-1}(X^x_s),U^x_s )-\lambda] \, ds - \int_t^T  Z^x_s \, dW^1_s-\int_t^T  U^x_s \, dW^2_s, \qquad   0 \leq t \leq T < \infty, 
\end{equation}
where, we recall,   $\lambda $ is a real number and it is part of the unknowns, and the equation has to hold for every $t$ and every $T$, see for instance \cite[section 4]{FuhHuTess}.
On the function $\psi: H \times \Xi^* \times H^* \to \mathbb{R}$ we assume:

\begin{itemize}
\item[$ {\bf (B.1)}$] $ (z,u) \to \psi(x,z,u)$ is a concave function at every  fixed $ x \in H$. 

Moreover there exist $L_x,L_z, L_u> 0$  such that 
\begin{equation}\label{Lip}
|\psi(x,z,u) - \psi(x',z'.u')| \leq L_x |x-x'| + L_z |z-z'|+L_u |u-u'|,\; \quad x,x'\in H, \ z,z'\in \Xi^*, \ u,u'\in H^*.
\end{equation}
Moreover $\psi(\cdot,0.0)$ is bounded. We denote $\sup_{x}|\psi(x,0.0)|$ by $M_\psi$.
\end{itemize}


We associate to $\psi$ its Legendre transformation (modified according to the fact that we are dealing with concave functions):
\begin{equation}
\psi^* (x,p,q)= \inf_{z \in \Xi^*, u\in H^*} \{ - z p -uq- \psi(x,z,u) \}, \qquad x \in H, p \in \Xi, q\in H.
\end{equation}

Clearly $\psi^*$ is concave w.r.t to $(p,q)$.

We collect some other properties of $\psi$ and $\psi^*$ we will use in the future:
\begin{Proposition}
Under hypothesis ${\bf (B.1)}$ we have that 
\begin{equation*}
\psi(x,z,u)= \inf_{(p,q) \in \mathcal{D}^*(x)}\{ -z p- uq- \psi ^* (x,p,q)\}.
\end{equation*}
where $ \mathcal{D}^*(x)= \{ (p,q): \psi^* (x,p,q) \not= -\infty \} \subset \left\{ (p,q)\in  \Xi\times H: |p|\leq L_z,\, |q|\leq L_u\right\}.$

Moreover  $ \mathcal{D}^*(x) = \mathcal{D}^*$ does not depend on $x \in H$, and there exists a $L_x >0$ such that
\begin{equation}\label{LipFenchel}
|\psi ^* (x,p,q)- \psi^*(x',p,q)| \leq  L_x |x-x'|, \qquad  x, x \in H, \; (p,q) \in  \mathcal{D}^*.
\end{equation}

Finally we remark that the above implies that for every $x \in H, z \in \Xi^*, u\in H^*$ :
\begin{equation*}  
\sup_{(p,q) \in  \mathcal{D}} \{  \psi(x,z,u) +   z p + uq + \psi^* (x,p,q) \}=0.
\end{equation*}
\end{Proposition}
\textbf{Proof.}
Since $\psi(x,\,\cdot\,,\,\cdot\,)$ is concave its double Legendre transform coincides with the function itself and the first relation follows immediately (see \cite{AubinMR549483}).

Then,  by the definition of $\psi^*$:
\begin{align*}
|\psi^* (x,p,q)-\psi^* (x',p,q) | \leq \sup_{z \in \Xi^*, \, u\in H^*}\left|  -zp -uq-\psi(x,z,u) +  zp +uq+\psi(x',z,u)\right| \leq L_{x}|x-x'|,
\end{align*}
thus we deduce that $\mathcal{D}^*$ doesn't depend on $x\in H$ and \eqref{LipFenchel} holds.
\ep

\medskip

As in \cite{FuhHuTess} we introduce, for each $\alpha >0$, the infinite horizon equation:
\begin{equation}\label{Yalpha}
Y^{x,\alpha}_t= Y^{x,\alpha}_T + \int_t^T [\psi(X^{x}_s, Z^{x,\alpha}_sG^{-1}(X^{x}_s),U^{x,\alpha}_s)-\alpha  Y^{x,\alpha}_s] \, ds - \int_t^T  Z^{x,\alpha}_s \, dW^1_s - \int_t^T  U^{x,\alpha}_s \, dW^2_s, 
\end{equation}
where $0 \leq t \leq T < \infty$.

The next result was proved in  \cite[Theorem 2.1]{Royer2004} when the $W$ is finite dimensional, the extension to the infinite dimensional case is straightforward, see also \cite[Lemma 4.2]{FuhHuTess}.
Notice that the random function, $ \widehat{\psi}(t,z,u) :=  {\psi}(X_t,G^{-1}(X_t) z,u)$, inherits the following properties:
\begin{equation}
|\widehat{\psi} (t,0,0)|=  | {\psi}(X_t,0,0)| \leq  M_{\psi},    \quad t \geq 0, \;  \mathbb{P}\hbox{- a.s.}.
\end{equation}
\begin{equation}
|\widehat{\psi} (t,z,u)- \widehat{\psi} (t,z',u')|  \leq  L_z M_{G^
{-1}}  |z-z'|+L_u|u-u'|\quad t \geq 0, \,  \;  \;z, z' \in \Xi^*, 
 \; u, u' \in H^*\ .
\end{equation} therefore it satisfies the assumptions in \cite[Lemma 4.2]{FuhHuTess}.
\begin{Theorem}
Let us assume $({\bf A.1}--{\bf A.5})$ and ${\bf (B.1)}$.
Then for every $\alpha >0$ there exists  a unique solution $( Y^{x,\alpha}, Z^{x, \alpha}, U^{x, \alpha} )$ to the BSDE \eqref{Yalpha} such that  $Y^{x,\alpha}$ is a bounded continuous process, $Z^{\alpha,x} \in L^{2, loc}_{\mathcal{P}}(\Omega\times [0,+\infty[;\Xi^*)$ and $U^{\alpha,x} \in L^{2, loc}_{\mathcal{P}}(\Omega\times [0,+\infty[;H^*)$.

Moreover  
\begin{equation}\label{boundY}
| Y^{x,\alpha}_t| \leq \frac{M_{\psi}}{\alpha}, \ \mathbb{P}  \text{-a.s.,  for all } t \geq 0.
\end{equation}
and
 \begin{equation}\label{boundZa}
\E \int_0^\infty   |e^{-\alpha s}Z^{x,\alpha}_s|^2 \, ds+ \E \int_0^\infty   |e^{-\alpha s}U^{x,\alpha}_s|^2 \, ds < \infty 
\end{equation}
\end{Theorem}

We define
\begin{equation}\label{valpha}
v^\alpha (x) = Y^{\alpha, x}_0
\end{equation}
The following is the main estimate of the paper.
\begin{Proposition}\label{main}
Under $({\bf A.1}--{\bf A.6})$ and ${\bf (B.1)}$ one has that for any $\alpha >0$:
\begin{equation}\label{Lipalpha}
|v^\alpha (x)- v^\alpha (x')| \leq  \frac{C}{\mu} |x-x'|, \qquad x,x'\in H.
\end{equation}
where $C$ depends on the constants in  $({\bf A.1}--{\bf A.5})$ and ${\bf (B.1)}$ but not on $\alpha$ (nor on $\mu$).
\end{Proposition}
\pr  \
Since, instead of  the pathwise  decay estimate holding  for $|X^x_t-X^{x'}_t|$ in the additive noise case (see \cite[Theorem 3.2]{FuhHuTess}), only the mean  bound \eqref{stimadiffSol} is true here we cannot proceed as in \cite[Theorem 4.4]{FuhHuTess}. Moreover, being the diffusion $X$, in general, degenerate, it is not possible to rely on the smoothing properties of its Kolmogorov semigroup (see \cite{Madec20151821}). On the contrary, concavity  assumption {\bf (B.1)} allows us to use control theoretic arguments.
\smallskip

First we notice that 
\begin{align*}
Y^{x,\alpha}_0= e^{-\alpha t} Y^{x,\alpha}_t + \int_0^t  e^{-\alpha s} \psi(X^{x}_s, Z^{x,\alpha}_sG^{-1}(X^{x}_s),U^{x,\alpha} _s)\, ds - \int_0^t  e^{-\alpha s}Z^{x,\alpha}_s \, dW^1_s- \int_0^t  e^{-\alpha s}U^{x,\alpha}_s \, dW^2_s
\end{align*}

Thus we have, taking also into account \eqref{boundY} and \eqref{boundZa}, that
\begin{equation}
Y^{x,\alpha}_0=   \int_0^ {+\infty} e^{-\alpha s} \psi(X^{x}_s, Z^{x,\alpha}_sG^{-1}(X^{x}_s),U^{x,\alpha}_s) \, ds 
- \int_0^{+\infty}  e^{-\alpha s}Z^{x,\alpha}_s \, dW^1_s
- \int_0^{+\infty}  e^{-\alpha s}U^{x,\alpha}_s \, dW^2_s.
\end{equation}
Moreover being $Y^{x,\alpha}_0$ deterministic, the  uniqueness in law for the system formed by equations \eqref{State1null} -\eqref{Yalpha} yields that it doesn't depend on the specific  independent Wiener processes.

We fix any stochastic setting $(\hat{\Omega},\hat{\mathcal{E}},(\hat{\mathcal{F}_t}), \hat{\mathbb{P}},(\hat{W_t}^1),(\hat{W_t}^2))$
where $((\hat{W_t}^1),(\hat{W_t}^2))$ are independent $(\hat{\mathcal{F}_t})$ Wiener processes with values in $\Xi$ and $H$ respectively.

Given any $(\hat{\mathcal{F}_t})$ progressively measurable process $\mathfrak{p}:=(p_t,q_t)$ with values in $\mathcal{D}^*$ by $(\hat{X}^{x,\mathfrak{p}}_t)$ we denote the unique mild solution  of the forward equation: 
\begin{equation}\label{State1bis}
d \hX^{x,\mathfrak{p}}_t \ = \ A \hX^{x,\mathfrak{p}}_tdt + F( \hX^{x,\mathfrak{p}}_t)dt+  Dq_tdt+Q G(\hX^{x,\mathfrak{p}}_t)p_tdt  + Q G(\hX^{x,\mathfrak{p}}_t)d \hat{W}^1_t + D  d \hat{W}^2_t  \qquad  \hat{X}^{x,\mathfrak{p}}_0  =  x.
\end{equation}
Clearly  $(\hat{X}^{x,\mathfrak{p}}_t)$ is also the unique mild solution of the forward equation: 
\begin{equation}\label{State2}
d \hX^{x,\mathfrak{p}}_t \ = \ A \hX^{x,\mathfrak{p}}_tdt + F( \hX^{x,\mathfrak{p}}_t)dt + Q G(\hat{X}^{x,\mathfrak{p}}_t)d \hat{W}^{1,\mathfrak{p}}_t + D  d \hat{W}^{2,\mathfrak{p}}_t  \qquad  \hat{X}^{x,\mathfrak{p}}_0 \ = \ x.
\end{equation}
where 
\begin{equation}
\hat{W}^{1,\mathfrak{p}}_t:= \hat{W}_t^1 +\int_0 ^t G^{-1}(\hX^{x,\mathfrak{p}}_s) p _s  \, ds ,\quad \hat{W}^{2,\mathfrak{p}}_t:= \hat{W}_t^2 +\int_0 ^tq _s  \, ds,
\end{equation}
and we know that under a suitable probability $\hat{\P}^\mathfrak{p}$ the processes $((\hat{W_t}^{1,\mathfrak{p}}),(\hat{W_t}^{2,\mathfrak{p}}))$ are independent  Wiener processes with values in $\Xi$ and $H$ respectively.\\
Let now  $(\hY^{x,\alpha,\mathfrak{p}}, \hZ^{x,\alpha,\mathfrak{p}},\hU^{x,\alpha,\mathfrak{p}} )$  be the solution to:
\begin{align*}
\hY^{x,\alpha,\mathfrak{p}}_t = &\hY^{x,\alpha,\mathfrak{p}}_T + \int_t^T \!\![\psi(\hX^{x,\mathfrak{p}}_s, \hZ^{x,\alpha,\mathfrak{p}}_sG^{-1}(\hX^{x,p}_s),\hU^{x,\alpha,\mathfrak{p}}_s)-\alpha  Y^{x,\alpha,\mathfrak{p}}_s] \, ds  \\ & - \int_t^T  \hZ^{x,\alpha,\mathfrak{p}}_s \, d\hat{W}^{1,\mathfrak{p}}_s- \int_t^T  \hU^{x,\alpha,\mathfrak{p}}_s \, d\hat{W}^{2,\mathfrak{p}}_s
\end{align*}
where  $0 \leq t \leq T < \infty $.

\noindent
By previous considerations one has, recalling that  
$
 \{  \psi(x,z) +   z p + uq + \psi^* (x,p) \} \leq 0, \forall  x \in H, z \in \Xi^*, u\in H^*, (p,q) \in \mathcal{D}^*,
$
that  for every $x \in H$
\begin{align*}
Y^{x,\alpha}_0 \!=& \hY^{x,\alpha,\mathfrak{p}}_0= 
 \\ =&  \int_0^{\infty}  \!\!\!\! e^{-\alpha s} \left[\psi(\hX^{x,\mathfrak{p}}_s, \hZ^{x,\alpha,\mathfrak{p}}_sG^{-1}(\hX^{x,\mathfrak{p}}_s), \hU^{x,\alpha,\mathfrak{p}}_s) +\hZ^{x,\alpha,\mathfrak{p}}_s G^{-1}(\hX^{x,\mathfrak{p}}_s) p _s+\hU^{x,\alpha,\mathfrak{p}}_s q_s  + \psi^* (\hX^{x,\mathfrak{p}}_s, p _s) \right] \! ds  \nonumber\\& -  \int_0^{+\infty}  e^{-\alpha s}\hZ^{x,\alpha,\mathfrak{p}}_s \, d\hat{W}^{1,\mathfrak{p}}_s  -  \int_0^{+\infty}  e^{-\alpha s}\hU^{x,\alpha,\mathfrak{p}}_s \,  d\hat{W}^{2,\mathfrak{p}}_s  -
 \int_0^{\infty}  \psi^* (\hX^{x,\mathfrak{p}}_s, p_s,q_s) \, ds  \\ \leq &  -  \int_0^{+\infty}  e^{-\alpha s}\hZ^{x,\alpha,\mathfrak{p}}_s \, d\hat{W}^{1\mathfrak{p}}_s  -  \int_0^{+\infty}  e^{-\alpha s}\hU^{x,\alpha,\mathfrak{p}}_s \,  d\hat{W}^{2,\mathfrak{p}}_s  -
 \int_0^{\infty}  \psi^* (\hX^{x,\mathfrak{p}}_s, p_s,q_s) \, ds .
\end{align*}
So:
 \begin{align}\label{stimaminorebis}
Y^{x,\alpha}_0  \leq  - \hat{\E}^{\mathfrak{p}} \int_0^{\infty}    e^{-\alpha s}\psi^*
 (\hX^{x,\mathfrak{p}}_s, p_s,q_s) \, ds
\end{align}
for arbitrary stochastic setting and arbitrary progressively measurable $\mathcal{D}^*$ valued control $\mathfrak{p}=(p,q)$.

\medskip
\noindent Then we fix $x\in H$ and assume, for the moment, that  $\forall \e \!>\!0$  there exists a stochastic setting

$$(\hat{\Omega}^{\e,x},\hat{\mathcal{E}}^{\e,x},(\hat{\mathcal{F}}^{\e,x}_t), \hat{\mathbb{P}}^{\e,x},(\hat{W}_t^{1,{\e,x}}),(\hat{W_t}^{2,\e,x}))$$ and a couple of  predictable processes $\mathfrak{p}^{\e,x}=(p^{\e,x}, q^{\e,x})$ with values in $\mathcal{D}^*$ such that (with the notations introduced above) the following holds $\mathbb{P}$ - a.s. for a.e. $s \geq 0$:
\begin{multline}\label{condepsottima}
\psi(\hX^{x,\mathfrak{p}^{\e}}_s,\hat{Z}^{x,\alpha,\mathfrak{p}^{\e,x}}_sG^{-1}(\hX^{x,\mathfrak{p}^{\e,x}}_s),\hat{U}^{x,\alpha,\mathfrak{p}^{\e,x}}_s) +\hat{Z}^{x,\alpha,\mathfrak{p}^{\e,x}}_s G^{-1}(\hX^{x,\mathfrak{p}^{\e,x}}_s) p^\e _s +\hat{U}^{x,\alpha,\mathfrak{p}^{\e,x}}_s q^{\e,x}_s \\ + \psi^* (\hat{X}^{x,\mathfrak{p}^{\e,x}}_s, p^{\e,x} _s, q^{\e,x} _s)
\geq - \e
\end{multline}

Proceeding as before we get:
\begin{align}\label{stimamaggiore}
Y^{x,\alpha}_0 =& \hY^{x,\alpha,\mathfrak{p}^{\e,x}}_0= \\ =&  \int_0^{\infty}  e^{-\alpha s} \left[\psi(\hX^{x,\mathfrak{p}^{\e,x}}_s, \hZ^{x,\alpha,\mathfrak{p}^{\e,x}}_sG^{-1}(\hX^{x,\mathfrak{p}^{\e,x}}_s), \hU^{x,\alpha,\mathfrak{p}^{\e,x}}_s) \right. \nonumber \\ &  \nonumber\qquad \qquad \left. +\hZ^{x,\alpha,\mathfrak{p}^{\e,x}}_s G^{-1}(\hat{X}^{x,\mathfrak{p}^{\e,x}}_s) p^{\e,x} _s+\hU^{x,\alpha,\mathfrak{p}^{\e,x}}_s q^{\e,x}_s  + \psi^* (\hX^{x,\mathfrak{p}^{\e,x}}_s, p^{\e,x} _s,q^{\e,x}_s) \right]\, ds  \nonumber\\&  \nonumber-  \int_0^{+\infty}  e^{-\alpha s}\hZ^{x,\alpha,\mathfrak{p}^{\e,x}}_s \, d\hat{W}^{1,\e,x}_s  -  \int_0^{+\infty}  e^{-\alpha s}\hU^{x,\alpha,\mathfrak{p}^{\e,x}}_s \,  d\hat{W}^{2,\e}_s  -
 \int_0^{\infty}  \psi^* (\hat{X}^{x,\mathfrak{p}^{\e,x}}_s, p^{\e,x}_s,q^\e _s) \, ds  \\ \geq &  -\frac{\e}{\alpha} -  \int_0^{+\infty}  e^{-\alpha s}\hZ^{x,\alpha,\mathfrak{p}^{\e,x}}_s \, d\hat{W}^{1,\e}_s  -  \int_0^{+\infty}  e^{-\alpha s}\hU^{x,\alpha,\mathfrak{p}^{\e,x}}_s \,  d\hat{W}^{2,\e,x}_s  -
 \int_0^{\infty}  \psi^* (\hat{X}^{x,\mathfrak{p}^{\e,x}}_s, p^{\e}_s,q^{\e,x} _s) \, ds\nonumber 
\end{align}
Thus  by  \eqref{stimaminorebis} taking into account \eqref{stimamaggiore} and  \eqref{LipFenchel}   we have:
\begin{align*}
Y^{x',\alpha}_0 -  Y^{x,\alpha}_0   \leq & \int_0^{\infty}    e^{-\alpha s} \hat{\E}^{\mathfrak{p}^{\e,x}}|\psi^* (\hat{X}^{x,\mathfrak{p}^{\e,x}}_s,  p^{\e,x}_s, q^{\e,x}_s)  -\psi^* (\hat{X}^{x',\mathfrak{p}^{\e,x}}_s, p^{\e,x}_s, q^{\e,x}_s) | \, ds + { \e}\\
&\leq \int_0^{\infty}    e^{-\alpha s}\,  \hat{\E}^{\mathfrak{p}^{\e,x}}|\hat{X}^{x,\mathfrak{p}^{\e,x}}_s -\hat{X}^{x',\mathfrak{p}^{\e,x}}_s | \, ds +\frac { \e}{\alpha},
\end{align*}
we stress the fact that  we keep the  stochastic setting $(\hat{\Omega}^{\e,x},\hat{\mathcal{E}}^{\e,x},(\hat{\mathcal{F}}^{\e,x}_t), \hat{\mathbb{P}}^{\e,x},(\hat{W}_t^{1,\e,x}),(\hat{W_t}^{2,\e,x}))$ and control $\mathfrak{p}^{\e, x}$ corresponding to the initial datum $x$ and just replace the initial state $x$ with a different one  $x'$.

Noticing now that both $(\hat{X}^{x,\mathfrak{p}^{\e,x}})$ and $(\hat{X}^{x',\mathfrak{p}^{\e,x}})$ satisfy (only the initial conditions differ):
$$d \hX_t \ = \ A \hX_tdt + F( \hX_t)dt+  Dq^{\e,x}_tdt+Q G(\hX_t)p^{\e,x}_tdt  + Q G(\hX^{x,\mathfrak{p}}_t)d \hat{W}^{1,{\e,x}}_t + D  d \hat{W}^{2,{\e,x}}_t$$
and taking into account \eqref{stimadiffSol} we can conclude that:
$$
Y^{x',\alpha}_0 -  Y^{x,\alpha}_0   \leq  L_x  \int_0^{\infty}   e^{-(\alpha + \frac{\mu}{2}) s} |x-x'| \, ds  +  \frac{ \e}{\alpha}\leq  \frac{C}{\mu} |x-x'| + \frac{ \e}{\alpha}.
$$

Interchanging the role of $x$ with $x'$
one gets:
\begin{align}\label{stimadifffinela1}
\left|Y^{x,\alpha}_0 -  Y^{x',\alpha}_0 \right|  &\leq \frac{C}{\mu} |x-x'| +\frac{ \e}{\alpha}.
\end{align}
where the constant $C$ is independent of $\alpha$, $\mu$ and $\e$  and is able to conclude \eqref{Lipalpha}  being $\e>0$ arbitrary.

\smallskip

We are left with the construction, for any fixed $x\in H$ and $\e>0$ 
of a stochastic setting $(\hat{\Omega}^{\e,x},\hat{\mathcal{E}}^{\e,x},(\hat{\mathcal{F}}^{\e,x}_t), \hat{\mathbb{P}}^{\e,x},(\hat{W}_t^{1,\e,x}),(\hat{W_t}^{2,\e,x}))$ and control $\mathfrak{p}^{\e,x}$ for which 
(\ref{condepsottima}) holds. 

We start from an arbitrary stochastic setting:
 $({\Omega},\mathcal{E},({\mathcal{F}}_t), {\mathbb{P}},({W}_t^{1}),({W_t}^{2}))$. Let $(X^x)$ be  the corresponding mild solution  of equation \eqref{State1} and $(Y^{x,\alpha},Z^{x,\alpha},U^{x,\alpha})$ the solution
 of \eqref{Yalpha}.
By a mea\-su\-rable selection argument see \cite[Theorem 4]{McSWar} we can find a couple of progressive measurable process $ \mathfrak{p}^{\e,x}=(p^{\e,x},q^{\e,x})$, (possibly depending on $\alpha$ as well), such that:
$$
\psi(X^{x}_s, Z^{x,\alpha}_sG^{-1}(X^{x}_s),U^{x,\alpha}_s) +Z^{x,\alpha}_s G^{-1}(X^{x}_s) p^{\e,x} _s +U^{x,\alpha}_s q^{\e,x}_s+ \psi_* (X^{x}_s, p^{\e,x} _s, q^{\e,x} _s)
\geq - \e.
$$
Then it is enough to set: 
\begin{equation}
\hat{W}^{1,\e,x}_t:= {W}_t^1 -\int_0 ^t G^{-1}(X^{x}_s) p^{\e,x} _s  \, ds ,\quad \hat{W}^{2,\e,x}_t:= {W}_t^2 -\int_0 ^tq^{\e,x} _s  \, ds,
\end{equation}
and choose $\hat{\Omega}^{\e,x}=\Omega$, $\hat{\mathcal{E}}^{\e,x}=\mathcal{E}$, $(\hat{\mathcal{F}}^{\e,x}_t))=(\mathcal{F}_t)$ and as
$\hat{\mathbb{P}}^{\e,x}$ the (unique) probability measure under which $((\hat{W}^{1,\e,x}_t), (\hat{W}^{2,\e,x}_t))$ are independent Wiener processes.
The claim then follows selecting the above control $\mathfrak{p}^{\e,x}$ and noticing that, by construction, $(\hat{X}^{x,\mathfrak{p}^{\e,x}})=(X^x)$. \ep

\medskip

Following \cite{FuhHuTess} we can find a function $ \bar{v}$ and a number  ${{\lambda}}$ such that:
\begin{equation} \label{defvbar}
 [v ^{\alpha_m}(x)- v ^{\alpha_m}(0) ] \to \bar{v}(x), \qquad  \forall x \in H,
 \end{equation}
 \begin{equation}\label{deflambda}
 \alpha_n v^{\alpha_m}(0) \to {\lambda}.
\end{equation}
where  $\{ \alpha_m\}_{m \in \mathbb{N}}$ is a suitable  subsequence constructed using a diagonal method.

We can then proceed as in \cite{FuhHuTess} to deduce from above the existence of a solution to \eqref{ergodicbsde} and the uniqueness of $\lambda$.

\begin{Theorem}\label{mainTheo}
Assume $({\bf A.1})-({\bf A.6})$ and ${\bf (B.1)}$, let ${\lambda}$ the number defined in \eqref{deflambda} and set 
$\bar{Y}^{x}_t := \bar{v}(X^x_t)$, where $\bar{v}$ is defined in \eqref{defvbar}. Then there exists $  \bar{Z}^x$ in $ L^{2, loc}_{\mathcal{P}}(\Omega\times [0,+\infty[;\Xi^*) $ and $  \bar{U}^x$  in  $L^{2, loc}_{\mathcal{P}}(\Omega\times [0,+\infty[;H^*) $ such that $(\bar{Y}^x,\bar{Z} ^x, \bar{U}^x, {\lambda})$ solves  equation \eqref {ergodicbsde},
 $\mathbb{P}$ -a.s. for all $ 0 \leq t \leq T$. 
 
 Moreover suppose that another quadruple $(Y' , Z',U', \lambda)$ where $Y'$  is a progressively measurable continuous
process verifying $|Y'_t|\leq c(1+|X^x_t|)$, ${Z}'  \in  L^{2, loc}_{\mathcal{P}}(\Omega\times [0,+\infty[;\Xi^*)$ , $ {U}' \in   L^{2, loc}_{\mathcal{P}}(\Omega\times [0,+\infty[;H^*) $ 
and $\lambda' \in \mathbb{R}$, 
satisfies \eqref{ergodicbsde}.
Then $\lambda' = \lambda$.

 Finally there exists a measurable function $\bar{\zeta}:  H \to \Xi^* \times H^* $ such that $(\bar{Z}^x_t, \bar{U}^x _t)=  \bar{\zeta}(X^x_t)$.
\end{Theorem}
\pr

Once \eqref{Lipalpha}, \eqref{defvbar} and \eqref{deflambda}
 are obtained, the proof  as far the first two statements is concerned follows exactly as in \cite [Theorem 4.4]{FuhHuTess}.


To get the existence of a function $\bar{\zeta}$, we proceed in the following way.
 For arbitrary fixed $0\leq t\leq T$ let  $(\bar{Y}^{x,t,T}, \bar{Z}^{x,t,T},\bar{U}^{x,t,T})$ be the solution to:
 \begin{equation}\label{system-funzionaleaddittivo}\left \{ \begin{array}{l}
  d X^{t,x}_s \ = \ A X^{t,x}_sds + F( X^{t,x}_s)ds + Q G(X^{t,x}_s)d W^1_s + Dd W^2_s , \\  X^{t,x}_t \ = \ x,\\
 - d {Y}^{x,t,T}_s =  \widehat{\psi}(X^{x,t}_s, Z^{x,t,T}_s,U^{x,t,T}_s ) \, ds - Z^{x,t,T}_s \, dW^1_s  - U^{x,t,T}_s \, dW^2_s- {\lambda} \, ds \\
 Y^{x,t,T}_T= \bar{v}(X^{x,t}_T)
\end{array}\right.
\end{equation}
Then we clearly have that $( \bar{Y}^x, \bar{Z} ^x, \bar{U}^x)$, restricted on $[0,T]$,  coincide with $(\bar{Y}^{x,0,T}, \bar{Z}^{x,0,T},\bar{U}^{x,0,T})$, for all $T>0$.
By \cite[Prop. 3.2]{Fuh2003} we know that there exists a measurable function $ \zeta^T : [0,T] \times H \to \Xi^*\times H^*$, such that  $ (\bar{Z}^{x,t,T}_s, \bar{U}^{x,t,T}_s) = \zeta^T (s, X^{x,t}_s), s \in [t,T]$. Moreover, see also \cite[Remark 3.3]{Fuh2003}, the map 
$[0,T]\ni (\tau,x)\rightarrow \zeta^T (\tau,x)$ is characterized in terms of the laws of $( \int_{\tau }^{\tau + \frac{1}{n}} \bar{Z}^{\tau,x,T}_s \, ds, \int_{\tau }^{\tau + \frac{1}{n}} \bar{U}^{\tau,x,T}_s \, ds)$, $ n \in \mathbb{N}$.

The uniqueness in law of the solutions to the system \eqref{system-funzionaleaddittivo}
together with the fact that its coefficients are time autonomous, we get:
\[   \int_{\tau }^{\tau + \frac{1}{n}} \bar{Z}^{\tau,x,T}_s \, ds \sim  \int_{0 }^{\frac{1}{n}} \bar{Z}^{0,x,T-\tau}_s \, ds  \sim \int_{0}^{ \frac{1}{n}} \bar{Z}^{x}_s \, ds   \]
and
\[   \int_{\tau }^{\tau + \frac{1}{n}} \bar{U}^{\tau,x,T}_s \, ds \sim  \int_{0 }^{ \frac{1}{n}} \bar{U}^{0,x,T-\tau}_s \, ds  \sim \int_{0}^{ \frac{1}{n}} \bar{U}^{x}_s \, ds   \]
So far we've proved that $\zeta^T(\tau,\cdot )$  does not depend neither from $T$  nor from $\tau$, thus we can define $\zeta^T(\tau,\cdot )=: \bar{\zeta}(\cdot)$  
and observe that  $(\bar{Z}^x_t, \bar{U}^x_t)=(\bar{Z}^{x,0,T}_t, \bar{U}^{x,0,T}_t)= \zeta^T(t,X^{x,0}_t )
=\bar{ {\zeta}}(X^x_t). $
\ep

%
%

\section{Ergodic Hamilton-Jacobi-Bellman}\label{sec_HJB}
Here we show that whenever $\bar{v}$ is  differentiable then $(\bar{v},\lambda)$ solves, in a {\em mild} form, the following {\em  Ergodic} HJB equation (see \cite{fuhrman2012stochastic}):
\begin{multline}\label{HJB}
\frac{1}{2}(\tr [QG(x)G^* (x) Q \nabla^2 \bar{v}(x)] +  \tr [DD^* (x) Q \nabla^2 \bar{v}(x))  + \langle  A x + F(x) , \nabla \bar{v}(x)\rangle =\\ - \psi(x,\nabla \bar{v}(x) Q,\nabla \bar{v}(x) D) + \lambda 
\end{multline}
mmoreover $\lambda$ characterize the ergodic limit of the parabolic solutions.

\smallskip 

We start by introducing the transition semigroup $(P_t)_{t \geq 0}$ corresponding to the diffusion $X^x$, see equation \eqref{State1null}:
\begin{equation}\label{transition}
P_t [\phi](x):= \E \, \phi(X^{x}_t), \qquad  \phi : H \to \R \text{ measurable and bounded. }
\end{equation}
We give the following definition, see \cite[Section 6]{fuhrman2012stochastic}:

\begin{Definition}\label{def-mild}
A pair $(v,\lambda)$ is a mild solution to the HJB equation \eqref{HJB} if  $ v \in \mathcal{G}^1 (H,\R)$ with bounded derivative and, for all $0\leq t \leq  T$, $x \in H$ it holds:
\begin{equation}
v(x) = P_{T-t}[v](x) + \int_t^T (P_{s-t}[\psi(\cdot, \nabla v(\cdot)Q, \nabla v(\cdot)D)](x)-\lambda) \, ds. 
\end{equation}
\end{Definition}
We have the following result.

\begin{Theorem}\label{teoident}
Assume    ${\bf(A.1--A.6)}$, ${\bf (B.1)} $ and that $\bar{v}$ is of class $\mathcal{G}^1$. Then $(\bar{v},{\lambda})$, defined in \eqref{defvbar}  is a {\em mild} solution of the HJB equation 
\eqref{HJB}. On the other hand if $ (v',\lambda')$ is a  mild solution of \eqref{HJB} then setting $ Y^x_t:= v'(X^x_t)$, $ Z^x_t
= \nabla v'(X^x_t) Q G (X_t^x) $ and $ U^x_t =  \nabla v'(X^x_t) D$,  we obtain that  $(Y^x, Z^x, U^x, \lambda)$ is a solution to equation 
\eqref{ergodicbsde}.

Moreover if $(v',\lambda')$ is another solution with $v'$ Gateaux differentiable with linear growth then $\lambda=\lambda'$.

Eventually, let for every $T >0$, $v^T(\cdot,\cdot)$ be the unique mild solution of the parabolic HJB equation:
\begin{multline}\label{HJBT}
\partial _t v^T (t,x) + \frac{1}{2}[{\tr}[QG(x)G^* (x) Q \nabla^2v^T (t,x)] +  \tr [DD^* (x) Q \nabla^2 v^T (t,x))  + \langle  A x + F(x) , \nabla v^T (t,x)\rangle =\\ - \psi(x,\nabla v^T (t,x) Q,\nabla v^T (t,x) D) , \qquad  v^T (T,x)=0.
\end{multline}
Then 
\begin{equation} \label{limite}
\lim_{T \to \infty} \frac{v^T (t,x) }{ T} = {\lambda}.
\end{equation}
\end{Theorem}
\pr  $ $ The existence part follows from  \cite [Theorem 6.2]{FuTessAOP2002}, while the uniqueness of $\lambda$ in the class of solutions that are Gateaux differentiable with linear growth follows as \cite[Theorem 4.6]{fuhrman2012stochastic}.
The only thing to prove is \eqref{limite}.

We prove \eqref{limite} in the case $t=0$. The general case follows in the same way just by replacing the initial time $0$ with $t$ in the forward equation \eqref{State1}.

We have that setting $\bar{Y}^{T,x}_s = v^T(s, X_s^x) -{\lambda} (T-s)$, $s\in [0,T]$, then $\bar{Y}^{T,x}$ solves:
\begin{equation}\label{parabolicT}\left \{ \begin{array}{l}
  - d {Y} ^{T,x} _s =  {\psi}(X^{x}_s,  G^{-1} (X^{x}_s)Z ^{T,x}_s ,  U ^{T,x}_s  ) \, ds - Z^{T,x}_s \, dW^1_s   - U^{T,x} _s \, dW^2_s- {\lambda} \, ds \\
\quad  \   Y^{T,x}_T= 0
\end{array}\right.
\end{equation}

Set $ \tilde {Y}^{T,x}_t= \bar{Y}^x_t -  \bar{Y}^{T,x}_t$, for all $ t\in [0,T]$,  then $\tilde{Y}^{T,x}$  verifies:
\begin{equation}\label{parabolicTdiff}\left \{ \begin{array}{l}
  - d\tilde {Y} ^{T,x} _s =  [{\psi}(X^{x}_s,  G^{-1} (X^{x}_s)\bar{Z} ^{x}_s ,  \bar{U} ^{x}_s  ) -  {\psi}(X^{x,t}_s,  G^{-1} (X^{x}_s)Z ^{T,x}_s ,  U ^{T,x}_s  )] \, ds -  ( \bar{Z}^x_s-Z^{T,x}_s) \, dW^1_s    \\ \qquad\qquad -  (\bar {U} ^x _s-U^{T,x} _s) \, dW^2_s- {\lambda} \, ds \\
\quad  \   \tilde{Y}^{T,x}_T= \bar{v}(X^{x}_T)
\end{array}\right.
\end{equation}

We rewrite  \eqref{parabolicTdiff}  as:

\begin{equation}\label{parabolicTdiff1}\left \{ \begin{array}{l}
  - d\tilde {Y} ^{T,x} _s =  \gamma^1_t   ( \bar{Z}^x_s-Z^{T,x}_s) \, ds  +  \gamma^2_t   ( \bar{U}^x_s-U^{T,x}_s) \, ds -  ( \bar{Z}^x_s-Z^{T,x}_s) \, dW^1_s    \\ \qquad\qquad -  (\bar {U} ^x _s-U^{T,x} _s) \, dW^2_s- {\lambda} \, ds \\
\quad  \   \tilde{Y}^{T,x}_T= \bar{v}(X^{x}_T)
\end{array}\right.
\end{equation}
where 
\begin{equation}\label{gamma1}
\gamma^1_s =
\left \{ \begin{array}{ll}
\frac{{\psi}(X^{x}_s,  G^{-1} (X^{x}_s)\bar{Z} ^{x}_s ,  \bar{U} ^{x}_s  ) -  {\psi}(X^{x}_s,  G^{-1} (X^{x}_s)Z ^{T,x}_s,  \bar{U} ^{x}_s) }{ |\bar{Z}^x_s-Z^{T,x}_s|^2_{\Xi^*}} (\bar{Z}^x_s-Z^{T,x}_s)^*& \text{if } \bar{Z}^x_s\not=Z^{T,x}_s,\\ 
0 & \text{ elsewhere}.
\end{array}\right.
\end{equation}
and 
\begin{equation}\label{gamma2}
\gamma^2_s =
\left \{ \begin{array}{ll}
\frac{{\psi}(X^{x}_s,  G^{-1} (X^{x}_s)Z ^{T,x}_s  ,  \bar{U} ^{x}_s  ) -  {\psi}(X^{x}_s,  G^{-1} (X^{x}_s)Z ^{T,x}_s ,  U^{T,x} _s))}{| \bar{U}^x_s-U^{T,x}_s|^2_H} (\bar{U}^x_s-U^{T,x}_s)^* & \text{if } \bar{U}^x_s\not=U^{T,x}_s,\\ 
0 & \text{ elsewhere}.
\end{array}\right.
\end{equation}
Hence, by a Girsanov argument,  we get that 
\begin{equation}
\tilde {Y}^{T,x}_0 = \E^{\gamma^1, \gamma ^2} (\bar{v}(X^{x}_T) )
\end{equation}
where the probability measure $\P^{\gamma^1, \gamma ^2}$ is the one under which $ W^{\gamma^1, \gamma ^2}_t= (W^1_t - \int_0^t \gamma^1_s \, ds, W^2_t - \int_0^t \gamma^2_s \, ds)$ is a cylindrical Wiener process  in $\Xi \times H$ in $[0,T]$.
Therefore by \eqref{stimaSol} and having  $\bar{v}$ Lipschitz,  we get that
\begin{equation}
\tilde {Y}^{T,x}_0 = \E^{\gamma^1, \gamma ^2} (\bar{v}(X^{x}_T) ) \leq  \kappa_{\gamma_1,\gamma_2} (1+ |x|)
\end{equation}
for some constant $\kappa_{\gamma_1,\gamma_2}$ independent of $T$.
Thus, noticing that   $ \tilde {Y}^{T,x}_0= \bar{v}(x) - v^T(0,x) + \lambda T$  we get that:
\begin{equation}
\lim_{T \to \infty } \frac{v^T(0,x)}{ T}= \lim_{T \to \infty } \frac{ \bar{v}(x)}{ T} + \lambda = \lambda.
\end{equation}

\ep

\section{Differentiability with respect to initial data}\label{sub-diff}

In this section we wish to present sufficient conditions under which  the function $\bar{v}$ defined in the section  above is differentiable. 

Throughout the section we assume   the following:
\begin{itemize}
\item [${\bf (C.1)}$] $F$ is of class $\mathcal{G}^1(H,H)$ and $G$ is of class  $\mathcal{G}^1(H,L(\Xi,H))$ 
\end{itemize}


We start from a straightforward result in the non-degenerate case.
 \begin{Proposition} Beside ${\bf(A.1--A.6)}$, ${\bf(B.1)}$ and  ${\bf(C.1)}$  assume that  the operator $\mathcal{Q}:=(Q,D): \Xi\times H\rightarrow H$  admits a right inverse $\mathcal{Q}^{-1}$ then $\bar{v}$  belongs to 
 class $\mathcal{G}^1(H)$.
 \end{Proposition}
 \pr
$ $ We fix $T>0$ and notice that $(\bar{Y}, \bar{Z}, \bar{U},\lambda)$ satisfies (see \eqref{ergodicbsde} and the definition of $\bar Y_t$ in Theorem \ref{mainTheo}):
$$
Y^x_t= \bar v(X^x_T) + \int_t^T [\widehat{\psi}(X^x_s,  \bar{Z}^x_s,\bar{U}^x_s )-\lambda] \, ds - \int_t^T  \bar{Z}^x_s \, dW^1_s-\int_t^T  \bar{U}^x_s \, dW^2_s, \qquad   0 \leq t \leq T < \infty, 
$$
where, we recall $\widehat{\psi}(x,z,u )=\psi(x,zG^{-1}(x),u )$ is lipschitz with respect to $z$ and $u$.
Moreover the forward equation \eqref{State1null} solved by $X^x$ can be rewritten as
$$
d X^x_t \ = \ A X^x_tdt + F( X^x_t)dt + \tilde{\mathcal{Q}}(X^x_t)d \mathcal{W}_t \qquad  X^x_0 \ = \ x.$$
where $\mathcal{W}_t :=\begin{pmatrix}
W^1_t \\
W^2_t
\end{pmatrix}$ is a $\Xi\times H$ valued Wiener process and $\tilde{\mathcal{Q}}(x)=(QG(x),D)$.

Under the present assumptions $\tilde{\mathcal{Q}}(x)$ turns out to be invertible with bounded right inverse:
$$[\tilde{\mathcal{Q}}(x)]^{-1}=\begin{pmatrix}
G^{-1}(x)& 0\\
0& I
\end{pmatrix} \mathcal{Q} ^{-1}$$
It is then straight forward to verify that all the assumptions in \cite[Theorem 3.10] {FuhTess2002BismutElworthy} are satisfied and consequently $\bar v$ (that coincides with the map $x\rightarrow \bar Y_x$) is in class $\mathcal{G}^1$
 \ep
 
 \medskip
 
 When the noise in the diffusion can be degenerate the situation is less simple and we will need quantitative conditions on the coefficients (see, for instance, \cite{Richou20092945}). 
 
\noindent We will now work under the  {\em joint dissipative condition}  ${\bf (A.7)}$  that, taking into account differentiability of $F$ and $G$ becomes:
 \begin{equation}\label{dissdiff}
2\langle  A y+ \nabla_x F(x)y, y\rangle_H  + || Q \nabla_x G(x) y ||^2_{L_2(\Xi,H)}\ \leq \ -\mu |y|_H^2,
\quad \forall y\in D(A), \,\forall x \in H.
\end{equation}

\medskip

Under the above assumptions the following well known  differentiability result for the forward equation \eqref{State1} holds:
\begin{Lemma} Under  ${\bf(A.1--A.5)}$,   ${\bf(A.7)}$ and ${\bf(C.1)}$ the map $x \to X^x$ is G\^ateaux differentiable. Moreover, for every $h \in H$, the directional derivative process $\nabla_x X^x h$, solves, $\mathbb{P}-$ a.s., the equation
\begin{equation}\label{eqfordir}
\nabla_x X^x_t h \ = \ e^{tA} h  + \int_0^t e^{(t-s)A}\nabla_x F( X_s^{x})\nabla_x X^x_s h  \, ds +  \int_0^t e^{(t-s)A}Q\nabla_x G( X_s^{x}) \nabla_x X^x_s h \, d W_s, \qquad t\geq0,
\end{equation}
Moreover 
\begin{equation} \label{stimagradXh}
\E |\nabla_x X^x_t h |^2 \leq  e^{-\mu t} |h|^2
\end{equation}
\end{Lemma}
\pr  $ $ Our hypotheses imply the Hypotheses 3.1 of  \cite{FuTessAOP2002}, therefore we can apply \cite[Prop 3.3]{FuTessAOP2002}. The estimate \eqref{stimagradXh} follows applying the It\^o formula to $| \nabla_x X^x_t h |^2$ and arguing as in Proposition \ref{Prop-esunstate}.
\ep 

$ $

We will need the following additional assumption to state the last result
\begin{itemize}
\item [${\bf (C.2)}$] $G$ and $G^{-1}$ are of class  $\mathcal{G}^1(H,L(\Xi))$ and $ \psi$ is of class $\mathcal{G}^1(H \times \Xi^*, \R)$
\end{itemize}
We eventually have:
\begin{Theorem}\label{maindiff}
Assume that  ${\bf(A.1--A.5)}$,  ${\bf(A.7)}$ and ${\bf (B.1)}$  hold with $\mu >  2({L^2_z} M_{G^{-1}}^2 + L^2_u )$,  moreover we assume ${\bf(C.1)}$ and  ${\bf(C.2)}$. Then the function $\bar{v}$ defined in \eqref{defvbar} is of class $\mathcal{G}^1(H,\R)$.
\end{Theorem}
\pr The proof is detailed in the Appendix. \ep

\section{Application to optimal control}\label{sec_control}
Let $\Gamma$ be a separable metric space,  an admissible control $\gamma$ is any $\mathcal{F}_t$ - progressively measurable $\Gamma$-valued process. 
The cost corresponding to a given control is defined as follows.
Let $ R_1: \Gamma \to \Xi$, $  R_2: \Gamma \to  H $ and $L: H\times \Gamma \to \R$ measurable functions such that, for some constant $c>0$, for all $ x,x'\in H$ and $\gamma\in \Gamma$:

\noindent
$\textbf{(E.1)}\qquad
|R_1(\gamma)|\leq c,\quad  |R_2(\gamma)|\leq c \quad |L(x,\gamma)| \leq c, \quad   |L(x,\gamma) -L(x',\gamma) |  \leq c |x-x'|.$

$ $

Let for every $x\in H$  be $X^x$ the  solution to \eqref{State1null}, then 
for every $ T>0$ and every control $\gamma$ we consider the Girsanov density:
\[ 
\rho^\gamma_T = \exp\left(   \int_0^T\!\!\! G^{-1}(X^x_s)R_1(\gamma_s) d W^1_s +  \int_0^T \!\!\!R_2(\gamma_s) \, dW^2_s -\frac{1}{2} \int_0^T[ | G^{-1}(X^x_s)R_1(\gamma_s)|^2 _{\Xi} + | R_2(\gamma_s)|^2 _{H}]\, ds \right) 
\] 
and we introduce the following ergodic cost corresponding to $x$ and $\gamma$:
\begin{equation*}\label{costo}
J(x,\gamma)= \limsup_{t \to \infty} \frac{1}{T} \, \E ^{\gamma,T} \!\!\int_0^T  L(X^x_s,\gamma_s)\, ds,
\end{equation*} 
where $\E^{\gamma,T}$ is the expectation with respect to $\mathbb{P}^\gamma : =\rho^{\gamma}_T \mathbb{P}$.
Notice that with respect to $\mathbb{P}^\gamma$ the processes $$W^{1,\gamma}_t:=-\int_0^t\!\!\! G^{-1}(X^x_s)R_1(\gamma_s)ds + dW^1_s,\quad
W^{2,\gamma}_t:=-\int_0^t\!\!\! R_2(\gamma_s)ds + dW^2_s$$
are independent cylindrical Wiener processes and with respect to them $X^x$ verifies:
$$\left\{ \begin{array}{ll}
d X^x _t  = \ A X^x_tdt + F( X^x_t)dt+ QR_1(\gamma_s)ds+DR_2(\gamma_s)ds + Q G(X^x_t)d W_t^{1,\gamma} + D  d W_t^{2,\gamma},  & t \geq 0, \\
 X^x_0 \ = \ x, &
\end{array}\right.
$$
and this justifies the above (weak) formulation of the control problem.

We introduce the {\em usual} Hamiltonian:
\begin{equation}\label{hamiltoniana}
\psi(x,z,u) = \inf_{\gamma \in \Gamma} \{ L(x,\gamma) + z R_1(\gamma) +u R_2(\gamma) \},\qquad  x\in H, z\in \Xi^*, u \in H^*
\end{equation}
that by construction is a concave function and, under {\bf (E.1)}, fullfils assumption {\bf (B.1)}. 
The forward backward system associated to this problem, is the following:
\begin{equation}\label{forbac}\left\{ \begin{array}{ll}
d X^x _t  = \ A X^x_tdt + F( X^x_t)dt + Q G(X^x_t)d W_t^1 + D  d W_t^2,  & t \geq 0, \\
 X^x_0 \ = \ x, & \\
-d Y^x_t = [\psi(X^x_t,  Z^x_tG^{-1}(X^x_t), U^x_t)-\lambda] \, dt -  Z^x_t \, dW_t^1 -  U^x_t \, dW_t^2. & 
\end{array}\right.
\end{equation}
By Theorem \ref{mainTheo} under ${\bf (A.1--A.6)}$  and ${\bf (E.1)}$ for every $x\in H$ there exists a  solution:
\begin{equation}\label{sol}
(\bar{Y}^x, \bar{Z}^x,   \bar{U}^x, {\lambda})= (\bar{v}(X^x), \bar{\zeta}_1(X^x),  \bar{\zeta}_2(X^x), {\lambda}), 
\end{equation}
where $\bar Y$ is a progressive measurable continuous process, $\bar Z \in L^{2, loc}_{\mathcal{P}}
(\Omega\times [0,+\infty[;\Xi^*)$, $\bar U \in L^{2, loc}_{\mathcal{P}}
(\Omega\times [0,+\infty[;H^*)$, $\lambda \in \R$, $\bar v$ is Lipschitz and $\bar{\zeta}_1$,   $\bar{\zeta}_2$ are measurable.

\medskip

Once we have  solved the above ergodic BSDE  the  proof of the following result containing the synthesis of the optimal control for the ergodic cost is identical to the one of \cite[Theorem 7.1]{FuhHuTess}.

\begin{Theorem}\label{teorema-main-controllo}
Assume ${\bf (A.1-- A.6)}$ and ${\bf (E.1)}$
Then the following holds:
\begin{itemize}
\item [(i)] For arbitrary control $\gamma$ we have $J(x,\gamma) \geq \lambda$, and equality holds if and only if the following holds  $\mathbb{P}$- a.s. for a.e. $t \geq 0$:
$$
L(X^x_t,\gamma_t)+\bar{ \zeta}_1(X^x_t) G^{-1}(X^x_t)R_1(\gamma_t) +  \bar{ \zeta}_2(X^x_t) R_2(\gamma_t)=  \psi(X^x_t, \bar{ \zeta}_1(X^x_t)G^{-1}(X^x_t), \bar{ \zeta}_2 (X^x_t)).$$

\item [(ii)] If the infimum is attained in \eqref{hamiltoniana} and 
$\rho: \Xi^* \times H^* \to $ is any measurable function realizing the minimum  (that always exists by Filippov selection theorem, see \cite{McSWar}) then the control $\bar{\gamma}_t= \rho (X^x_t,\bar{ \zeta _1}(X^x_t ),  \bar{ \zeta _2}(X^x_t ))$ is optimal, that is $J(x,\bar{\gamma})=\lambda.$
\item [(iii)] Finally if $\bar{v}$ is in class $\mathcal{G}^1$ then it is a mild solution of equation \eqref{HJB} and  $ \bar{ \zeta}_1 
= \nabla \bar{v} Q G$ and  $ \bar{ \zeta}_2
= \nabla \bar{v} D$ .
\end{itemize}
\end{Theorem}

\subsection{Examples}\label{sub-ex}

\begin{Example}\label{esempio_bordo} { \rm We consider an ergodic control problem for a stochastic heat equation  controlled through the boundary
\begin{equation}\label{eq_esempiobordo}
\left\{  \begin{array}{ll}
d_t x(t, \xi)=  \frac{\partial}{\partial \xi ^2}x(t, \xi) \, dt  + d(\xi) \dot{\mathcal{W}}(t, \xi) \, dt, & t \geq 0, \ \xi \in (0,\pi), \\ 
x(t,0)=  y(t), \qquad x(t,\pi)= 0,  \\
x(0,\xi)= x_0(\xi),& \xi\in (0,\pi) \\

d y (t)= b(y(t))\, dt +  \sigma(y(t))\rho(\gamma(t))dt+\sigma(y(t)) \, d B_t, & t \geq 0,\\
y(0)= x \in \R.
\end{array}
\right.
\end{equation}
where $\mathcal{W}$ is the space-time white noise  on $[0,+\infty) \times [0,\pi]$ and $B$ is a brownian motion.  An admissible control $\gamma$ is a predictable process $ \gamma : \Omega \times [0,+\infty)  \to \R$.
The cost functional is
\begin{equation}\label{costobordo}
J(x_0,\gamma)= \liminf_{T \to +\infty} \, \frac{1}{T}\,  \E \int_0^T \int_0^\pi \ell(x(t, \xi), \gamma(t)) \, d \xi \, dt.
\end{equation}
We assume that
{\rm
\begin{enumerate}
\item $b : \R \to \R$ is a measurable function such that
\begin{align*}
|b(y)-b(y') | &\leq  L_b |y-y'|,
\end{align*}
for a suitable positive constant  $L_b$, for  every $ y,y \in \R$. 

\item  $ \sigma:  \R \to \R$ is a measurable and bounded function, such that 
\begin{align*}
|\sigma(y)-\sigma(y') | &\leq  L_\sigma  |y-y|,
\end{align*}
 for suitable positive constants  $L_\sigma$ and there exists a  suitable positive $ \delta $ such that:
\[ |\sigma (y))| \geq \delta > 0,
\]
for every $y \in \R$.
\item 
there exists $\mu>0$ such that for all $y,y' \in \R$:
\begin{equation}\label{condissesempio}
2\langle b(y)- b(y'), y - y' \rangle + | \sigma(y) - \sigma(y') |^2 \ \leq \ -\mu|y - y'|^2,
\end{equation}

\item $d: [0,\pi] \to \R $, $\rho: \R \to \R$ are bounded and measurable functions.
\item $\ell: \R \times \R \to \R$ is a measurable and bounded function such that
\[ |\ell(x, \gamma)- \ell(x',\gamma)| \leq L | x-x'|,
\]
 for a suitable positive constant $L$, for  every $x, x', \gamma \in \R$.
\end{enumerate}}

Under these hypotheses, see \cite{LasTri}, the above equation can be reformulated in an infinite dimensional space as:
\begin{equation}\label{eq_esempiobordo-secondo}
\left\{  \begin{array}{ll}
d_t \mathcal{X}_t=  \Delta \mathcal{X}_t \, dt -  \Delta \mathfrak{r} y(t)dt+ \tilde{D} d\tilde{W}_t \,, & t \geq 0, \ \xi \in [0,\pi], \\
\mathcal{X}_0= x_0(\cdot), & \xi\in (0,\pi)\\ 
d y (t)= b(y(t))\, dt +  \sigma(y(t))\rho(u(t))dt+\sigma(y(t)) \, d B(t), & t \geq 0,\\
y(0)= y_0\in \R.
\end{array}
\right.
\end{equation}
where $\mathcal{X}_t:=x(\cdot)$ is in $L^2(0,\pi)$, $\tilde{W }$ is a cylindrical Wiener process in $L^2(0,\pi)$,  $\tilde{D}$ is the bounded operator in $L^2(0,\pi)$ corresponding to multiplication by a bounded function $d$,
 $\Delta$ is the realisation of the Laplace operator  with Dirichlet boundary conditions in $L^2(0,\pi)$,  that is (denoting by $\mathcal{D}(\Delta)$ the domain of the operator)
 $$\mathcal{D}(\Delta) = H^2(0,\pi) \cap H^1_0(0,\pi),\quad \;\;\Delta f=\frac{\partial^2 f}{\partial \xi^2 },\;\; \forall f\!\in\! \mathcal{D}(\Delta) $$

  Finally $\r(\xi)= 1-\frac{\xi}{\pi}, \ \xi \in [0,\pi]$ is the 
solution to
\begin{equation}\left\{ \begin{array}{ll}
 \frac{\partial^2 \r}{\partial \xi^2 } (\xi)= 0, & \xi \in (0,\pi) ,\\ 
 \r(0)=1, \qquad   \r (\pi)=0.
\end{array}\right.
\end{equation}
It is well known that $\Delta$ generates an analytic semigroup of contractions (of negative type  $-1$) moreover, for any $ \delta >0$, $\mathfrak{r}\in \mathcal{D}((-\Delta)^{1/2-\delta})$ (where $(-\Delta)^{\alpha}$  denotes the fractional power). Standard results on analitic semigroups then yield:

\begin{equation}	\label{stimafrazionaria}
| (-\Delta ) e^{t \Delta } \r |_{L^{2}(0,\pi)}
 \leq c_{\r} e^{ -t}t ^{- (\frac{1}{2}+\delta)}, \qquad t > 0.
\end{equation}

\noindent We are now in a position to  rephrase the  problem according to our general framework. Indeed setting $H= L^2(0,\pi) \times \R$, $\Xi=\mathbb{R}$ and
 $X_t=\begin{pmatrix} \mathcal{X}_t,y(t)
\end{pmatrix}$ equation \eqref{eq_esempiobordo-secondo} becomes

\begin{equation}\left\{ \begin{array}{ll}
d X^x _t  = \ A X^x_tdt + F( X^x_t)dt + Q G(X^x_t)\rho(\gamma_t)dt + Q G(X^x_t)d W_t^1 + D  d W_t^2,  & t \geq 0, \\
 X^x_0 \ = \ x.
\end{array}\right.
\end{equation}
where:

\begin{enumerate}
\item $A = \begin{pmatrix}  -\Delta  & -\Delta R \\
0 & 0 \end{pmatrix}$ where $R: \mathbb{R} \to D((-\Delta) ^{\frac{1}{2}-\delta})$, is defined as $Ry= \r(\cdot)y$, $y\in \mathbb{R}$

It is easy to verify that $A$ generates a $C_0$-semigroup in $H$.

\item
$F: H \to H$, is defined as:
$F  \begin{pmatrix}\mathcal{X}\\y \end{pmatrix}= \begin{pmatrix} 0 \\ b(y)\end{pmatrix}$,

$Q: \Xi \to H$ is defined as:
$ Q  y   =  \begin{pmatrix}0\\ y \end{pmatrix}$, 

$G: \Xi \to \Xi$, is defined as:
$ G( y)  = \sigma(y)$

$D: H\to H$ is
 defined as:
$ D  \begin{pmatrix}\mathcal{X}\\y \end{pmatrix}=   \begin{pmatrix}  \tilde{D}\mathcal{X} \\0 \end{pmatrix}.$
\item $W^1(t)=B(t)$ and
 $( W^2) $ is a cylindrycal Wiener process in $H$.


\end{enumerate}

Hypotheses ${\bf (A.1--A.5) }$ are immediately verified, we have to check ${\bf (A.6)}$.  We come back to the formulation \eqref{eq_esempiobordo-secondo} and start  with  the second component $y$ (that only depends on $y_0$).
By \eqref{condissesempio},  Proposition \ref{propdiss} gives:
\begin{align}\label{stima_seconda_coordinata}
\E | y^{y_0}(t)- y^{y'_0}(t) | ^2 \leq e ^{-2 \mu t} |y_0-y_0'|^2.
\end{align}
Coming now to the first component we have that it fullfills in $L^2(0,\pi)$ the following mild formulation:
\begin{align*}
\mathcal{X}^{x_0,y_0}_t = e^{t\Delta } x_0- \int_0^t  \left[\Delta  e^{ (t-s) \Delta  } \r \right] y^{y_0} (s) \, ds +  \int_0^t  e^{ (t-s) \Delta  } D \, dW_s
\end{align*}
Thus considering two different initial data
\begin{align*}
\mathcal{X}^{x_0,y_0}_t-  \mathcal{X}^{x_0',y_0'}_t = e^{ t\Delta  } (x_0-x'_0)- \int_0^t \Delta  e^{ (t-s) \Delta  }( \r y^{y_0} (s)-  \r y^{y'_0} (s))  \, ds.
\end{align*}
By \eqref{stimafrazionaria} and \eqref{stima_seconda_coordinata} choosing $\mu_0\in (0,1 \wedge \mu)$
\begin{align*} 
 \E  |\mathcal{X}^{x_0,y_0}_t-  \mathcal{X}^{x_0',y_0'}_t| & \leq e^{ -t} |x_0-x_0'| + \int_0^t   e^{- (t-s)}(t-s) ^{- (\frac{1}{2}+\delta)}  e^{-\mu s}|y_0-y'_0| \, ds  \\
& \leq e^{ -t} |x_0-x_0'|  + e^{ -\mu_0 t}\left[ \int_0 ^t  e^{-(1  - \mu_0)(t-s)} (t-s) ^{- (\frac{1}{2}+\delta)}\, ds\right] |y_0-y'_0| .
\end{align*}

That implies that  \eqref{stimadiffSol} holds. In the same way one gets the proof of \eqref{stimaSol}.

We notice that it is  not at all obvious that the stronger versions \eqref{stimaSolqua}, \eqref{stimadiffSolqua} holds in this case.

As far as the control functional is concerned it is enough to set 
$L(X,\gamma)=\int_0^{\pi}\ell(\xi,\mathcal{X}(\xi),\gamma)d \, \xi$ and to verify in a straightforward way that \textbf{(E.1)} holds (in this case $R_1=\rho$, $R_2=0$, $\Gamma=\mathbb{R}$).

Thus all the hypotheses of Theorem \ref{teorema-main-controllo} hold and points (i) and (ii) in its thesis give the optimal ergodic cost and strategy in terms of the solution to the ergodic BSDE  in \eqref{forbac}.

}

\end{Example}
\begin{Example}\label{esempioaltro} { \rm We consider an ergodic control problem for a stochastic heat equation with Dirichlet boundary conditions with nonlinearity controlled through a one dimensional process $y$.
\begin{equation}\label{eq_esempidue}
\left\{  \begin{array}{ll}
d_t x(t, \xi)=  \frac{\partial}{\partial \xi ^2}x(t, \xi) \, dt  +f(x(t, \xi),y(t))+ d(\xi) \dot{\mathcal{W}}(t, \xi) \, dt, & t \geq 0, \ \xi \in (0,1), \\ 
x(t,0)= x(t,1)= 0,  \\
x(0,\xi)= x_0(\xi),& \xi\in (0,1) \\

d y (t)= b(y(t))\, dt +  \sigma(y(t))\gamma(t)dt+\sigma(y(t)) \, d B_t, & t \geq 0,\\
y(0)= y_0 \in [-1,1].\end{array}
\right.
\end{equation}
where $\mathcal{W}$ is the space-time white noise  on $[0,+\infty) \times [0,1]$ and $B$ is a brownian motion.  An admissible control $\gamma$ is a predictable process $ \gamma : \Omega \times [0,+\infty)  \to [-1,1]$.
The cost functional is
\begin{equation}\label{costoaltro}
J(x_0,\gamma)= \liminf_{T \to +\infty} \, \frac{1}{T}\,  \E \int_0^T \left[\int_0^1 (\ell( x(t, \xi),y(t))d\xi+\gamma^2(t)\right]\, dt.
\end{equation}
We assume:

\begin{enumerate}
\item $f : \R^2 \to \R$ is a Lipschitz map.
We fix two constants $L_f>0$ and $\mu_f\in \mathbb{R}$ such that
\begin{align*}
|f( x,y)-f( x',y) | &\leq  L_f  (|x-x'| + |y-y' |),\quad \langle f(x,y)-f(x,y'),  x-x' \rangle  &\leq-  \mu_f  |x-x'| ^2,
\end{align*}
for every $x, x', y, y \in \R$.

\item $b :  \R \to \R$ is Lipschitz. We fix a constant $\mu_b\in \mathbb{R}$ such that:
\begin{align*}
\langle b(y)-b(y'), y-y' \rangle  &\leq  - \mu_b  |y-y'| ^2, \qquad \forall  y,y' \in \R
\end{align*}
\item  $ \sigma:  \R^2 \to \R$ is a Lipschitz and bounded. We fix $L_\sigma $ such that 
\begin{align*}
|\sigma(y)-\sigma(y') | &\leq  L_\sigma |y-y'|, \qquad \forall  y,y'\in \R,
\end{align*}
We also assume that there exists a  suitable positive $ \delta $ such that:
\[ |\sigma (y))| \geq \delta > 0,  \qquad \forall  y\in \R.
\]
\item $d: [0,1] \to \R $ is a bounded and measurable function.
\item $\ell:  \R^2 \to \R$ is bounded and Lipschitz 

\end{enumerate}

As in the previous example  the above equation can be reformulated in an infinite dimensional space as:
$$
\left\{  \begin{array}{ll}
d_t \mathcal{X}_t=  \Delta \mathcal{X}_t \, dt + f(\mathcal{X}_t,y(t))dt+ \tilde{D} d\tilde{W}_t \,, & t \geq 0, \ \xi \in [0,1], \\
\mathcal{X}_0= x_0(\cdot), & \xi\in [0,1]\\ 
d y (t)= b(y(t))\, dt +  \sigma(y(t))\gamma(t)dt+\sigma(y(t)) \, d B(t), & t \geq 0,\\
y(0)= y_0\in \R.
\end{array}
\right.
$$
where $\mathcal{X}_t:=x(\cdot)$ is in $L^2(0,1)$, $\tilde{W }$ is a cylindrical Wiener process in $L^2(0,1)$,   $\Delta$ is the realisation of the Laplace operator  with Dirichlet boundary conditions in $L^2(0,1)$, $\tilde{D}$ is the bounded operator in $L^2(0,1)$ corresponding to multiplication by a bounded function $d$.

Finally setting $H= L^2(0,1) \times \R$, $\Xi=\mathbb{R}$,  $\Gamma = [-1,1]$  and
 $X_t=\begin{pmatrix} \mathcal{X}_t,y(t)
\end{pmatrix}$ equation \eqref{eq_esempiobordo} becomes

\begin{equation}\left\{ \begin{array}{ll}
d X^x _t  = \ A X^x_tdt + F( X^x_t)dt + Q G(X^x_t)\gamma_tdt + Q G(X^x_t)d W_t^1 + D  d W_t^2,  & t \geq 0, \\
 X^x_0 \ = \ x.
\end{array}\right.
\end{equation}
and the cost takes our general form:
$$
J(x_0,\gamma)= \liminf_{T \to +\infty} \, \frac{1}{T}\,  \E \int_0^T L(X(t),\gamma(t))\, dt.
$$
where
\begin{enumerate}
\item $A = \begin{pmatrix}  -\Delta  & 0 \\
0 & 0 \end{pmatrix}$ 
 generates a $C_0$-semigroup in $H$.
We also have that 
\[ \langle   A X, X\rangle_ H = \langle \Delta \mathcal{X}, \mathcal{X} \rangle _ {L^2(0,1)}\leq -\mu_{\Delta} |  \mathcal{X} |^2 _{L^2(0,1)}, \]
for some $\mu_{\Delta} >0$.
\item
$F: H \to H$, is defined as:
$F  \begin{pmatrix}\mathcal{X}\\y \end{pmatrix}= \begin{pmatrix}  f(\mathcal{X},y)\\ b(y)\end{pmatrix}$,

$Q: \Xi \to H$ is defined as:
$ Q  y   =  \begin{pmatrix}0\\ y \end{pmatrix}$, 

$G: \Xi \to \Xi$, is defined as:
$ G( y)  = \sigma(y)$

$D: H\to H$ is
 defined as:
$ D  \begin{pmatrix}\mathcal{X}\\y \end{pmatrix}=   \begin{pmatrix}  \tilde{D}\mathcal{X} \\0 \end{pmatrix}.$
\item $W^1(t)=B(t)$ and
 $( W^2) $ is a cylindrycal Wiener process in $H$.
\item $ L: H \times \Gamma \to \R$,   $L(X, \gamma)=\displaystyle \int_0^{1}\ell(\mathcal{X}(\xi),y)d \, \xi + |\gamma|^2$ 
\end{enumerate}
We
 also notice that in this case the Hamiltonian defined as in \eqref{hamiltoniana} becomes:
\begin{equation}\label{hamintoniana_esempio_due}
\psi\left(\begin{pmatrix}
\mathcal{X}\\y
\end{pmatrix},z\right)= -\frac{z^2}{4}I_{[-2,2]}(z)+(1-|z|)I_{[-2,2]^c}(z)+\int_0^{1}\ell(\mathcal{X}(\xi),y)d \, \xi
\end{equation}
We also assume that there exists $\bar{\mu}>0$  such that
\begin{equation}\label{definitanegativa}
\begin{pmatrix}
-\mu_\Delta -\mu_f &\frac{1}{2}L_f\\
\frac{1}{2}L_f&  -\mu_b +\frac{1}{2} {L_\sigma}
\end{pmatrix}\leq  -\bar{\mu} \,I_{\mathbb{R}^2}
\end{equation}

Hypotheses ${\bf (A.1--A.5) }$ are immediately verified. Moreover
relation \eqref{definitanegativa} ensures that  $\bf{(A.7)}$  holds as well. Finally \textbf{(E.1)} is straight forward (in this case $R_1= id$, $R_2=0$). Thus the hypotheses of Theorem \ref{teorema-main-controllo} hold and points (i) and (ii) in its thesis give the optimal ergodic cost and strategy in terms of the solution to the ergodic BSDE  in \eqref{forbac}.

We finally wish to apply  the differentiability 
result in Theorem \ref{maindiff} to this specific example. We notice that by \eqref{hamintoniana_esempio_due} the Hamiltonian $\psi$ is concave and differentiable with respect to $z$ with $\nabla_{z}\psi\leq 1$. Thus $\textbf{(B.1)}$ holds and we can choose $L_z=1$ in \eqref{Lip}.
 If we assume  that $ f$ $b$ $\sigma$ and $\ell$ are of class $C^1$ in all their variables then  \textbf{(C.1)} and  \textbf{(C.2)}  hold, moreover if we impose that $\bar{\mu}>2\delta^{-2}$ (here, comparing with Theorem \ref{maindiff}, $L_u=0$, $M_{G^{-1}}=\delta^{-1}$) then all the assumptions of Theorem \ref{maindiff} are verified and we can conclude that function $\bar{v}$ in  Theorem \ref{teorema-main-controllo} is differentiable. Consequently  point (iii) in Theorem \ref{teorema-main-controllo} as well applies here and we obtain that $\bar{v}$ is a mild solution of equation \eqref{HJB} and that the optimal feedback law can be characterized in terms of the gradient of $\bar{v}$.

}\end{Example}

\appendix
\section{Proof of Theorem \ref{maindiff}}

\smallskip

We will need to use  some results from \cite[Theorem 5.21 and Section 5.6]{ParRasBook}. The first concerns finite horizon BSDEs and the estimate of their solution, while  the  second  concerns the infinite horizon case. We restate them in  our setting as follows:
\begin{Lemma} \label{ParRasFin}
Let us consider the following equation:
\begin{align}\label{parras1}
-d\, Y_t  & = ( \phi(t, Z_t, U_t) \, dt  -\alpha  Y_t) \, dt -   Z_t \, dW_t^1-  U_t \, dW_t^2,  \qquad   Y_T = \eta , \qquad t \in [0,T],\  \alpha \geq 0.
\end{align}
assume that:

\begin{enumerate}
\item $|  \phi (t,z,u)- \phi (t,z',u')| 
\leq \ell (t) ( |z-z'|^2 + |u-u'| ^2)^{1/2}$,  $\forall z,z'\in \Xi ^* , u,u'\in H^*$, $\mathbb{P}-a.s.$
for some $\ell \in  L^2([0,T])$;
\item for $\nu_t:= \displaystyle\int_0 ^t \ell^2(s)\, ds$, one has
\begin{equation}\label{ipoVtfin}
\E \left( e^{2 \nu_T -2 \alpha T} |\eta|^2 \right) < \infty, \qquad  \E \left(\int_0^T e^{\nu_s - \alpha s} |\phi(s,0,0)| \, ds \right)^2 < \infty.
\end{equation}
\end{enumerate}
Then there exists a unique solution $(Y,Z,U) \in L^{2}_\Pc (\Omega;{C} ([0,T];\mathbb{R})) \times  
L^2_\Pc (\Omega \times [0,T];\Xi ^*)\times  
L^2_\Pc (\Omega \times [0,T];H ^*)$ and it verifies for all $0\leq t\leq T$:
\begin{align}\label{stimaparrasfin}
&\E^{\mathcal{F}_t} (\sup_{s \in [t,T]}  e^{2(\nu_s -\alpha s)} |Y_s |^2) + \E ^{\mathcal{F}_t}\left( \int_t^T e^{2(\nu_s -\alpha s)} |Z_s |^2\, ds\right) +  \E ^{\mathcal{F}_t}\left( \int_t^T e^{2(\nu_s -\alpha s)} |U_s |^2\, ds\right) \nonumber
\leq  \\&
\E^{\mathcal{F}_t} \left( e^{2 \nu_T - 2\alpha T} |\eta|^2 \right) + \E ^{\mathcal{F}_t}\left(\int_t^T e^{V_s - \alpha s} |\phi(s,0,0)| \,ds \right)^2, \qquad \mathbb{P}-a.s., \quad \   t \in [0,T]
\end{align}
\end{Lemma}
\begin{Lemma} \label{ParRasInf}
Let us consider the following equation for $\alpha \geq 0$:
\begin{align}\label{parras2}
-d\, Y_t  & = ( \phi(t, Z_t, U_t) \, dt  -\alpha  Y_t) \, dt -   Z_t \, dW_t^1 -  U_t \, dW_t^2, \qquad t \geq 0, \  .
\end{align}
Assume that:
\begin{enumerate}
\item $|  \phi (t,z,u)- \phi (t,z',u')| 
\leq \ell (t) ( |z-z'|^2 + |u-u'| ^2)^{1/2}$,  $\forall z,z'\in \Xi ^* , u,u'\in H^*$, $\mathbb{P}-a.s.$
for some $\ell \in  L^2_{loc}([0,+\infty [)$;
\item for $\nu_t:=  \displaystyle \int_0 ^t \ell^2(s)\, ds$, one has
\begin{equation}\label{ipoVtinf}
 \E \left(\int_0^{\infty} e^{\nu_s} |\phi(s,0,0)| \, ds\right)^2 < \infty.
\end{equation}
\end{enumerate}
Then there exists a unique triple of processes $(Y,Z, U )$ with 
 $ Y\in L^{2, loc}_\Pc (\Omega;{C} ([0,+\infty[;\mathbb{R}))$, $Z\in L^{2,loc}_\Pc (\Omega \times [0,+\infty [;\Xi ^*) $, $U\in
L^{2,loc}_\Pc (\Omega \times [0,+\infty [;H ^*)$, such that
\begin{equation}
\E (\sup_{t \in [0,T] } e^{2 \nu_t } |Y_t |^2) < + \infty,  \;\forall T \geq 0,\qquad 
\lim_{T\to \infty}\E (e^{ 2 \nu_T } |Y_T |^2) =0.
\end{equation}
Moreover
\begin{multline}\label{stimaparrasinf}
\E^{\mathcal{F}_t} (\sup_{s \geq t } e^{2 \nu_s } |Y_s |^2) + \E ^{\mathcal{F}_t}\left( \int_t^{\infty} e^{2\nu_s } (|Z_s |^2 + |U_s|^2)\, ds\right)     \leq  C \, 
\E^{\mathcal{F}_t} \left(\int_t^{\infty} e^{  \nu_s} |\phi(s,0,0)| \,ds \right)^2,  \,  \mathbb{P}-a.s.\end{multline}
for some  positive constant $C$.
\end{Lemma}


{\bf Proof of Theorem \ref{maindiff}}. The proof is split into two parts. The first deals with approximating functions  $v^{\alpha}$ defined in \eqref{valpha}

\medskip 
\noindent\textit{Part I - Differentiability of $v^{\alpha}$}

\noindent We first have to come back to the elliptic approximations:
\begin{equation}\label{Yalpha1}
Y^{x,\alpha}_t= Y^{x,\alpha}_T + \int_t^T [\psi(X^{x}_s, Z^{x,\alpha}_sG^{-1}(X^{x}_s), U^{x,\alpha}_s)-\alpha  Y^{x,\alpha}_s] \, ds - \int_t^T  Z^{x,\alpha}_s \, dW^1_s - \int_t^T  U^{x,\alpha}_s \, dW^2_s,   
\end{equation}
and for those equations we prove that:
\begin{Proposition}\label{propdiff}
Under the same assumptions of Theorem \ref{maindiff} we have that, for each $\alpha>0$,  the map $ x\to Y^{x,\alpha}_0$ belongs to $\mathcal{G}^1(H,\R)$.
\end{Proposition}
\pr  $ $ We fix $n\in \mathbb{N}$ and  introduce the following finite horizon approximations where  $0 \leq t \leq n   $:
$$
Y^{x,\alpha,n}_t\!=\!  \int_t^n [\psi(X^{x}_s, Z^{x,\alpha,n}_sG^{-1}(X^{x}_s), U^{x,\alpha,n}_s)-\alpha  Y^{x,\alpha,n}_s] \, ds - \int_t^n \! Z^{x,\alpha,n}_s \, dW^1_s - \int_t^n \! U^{x,\alpha,n}_s \, dW^2_s.
$$
For such equations \cite[Prop. 3.2]{HuTess} holds true, moreover we have from \cite[Propositions 5.6 and 5.7]{FuTessAOP2002} that $x \to Y^{x,\alpha,n}_0:= v^{\alpha, n}(x)$  belongs to $\mathcal{G}^1(H,\R)$ and  $Z^{x,\alpha,n}_t= \nabla _x v^{\alpha, n}(X_t^x) G(X^{x}_t)$ and  $U^{x,\alpha,n}_t= \nabla _x v^{\alpha, n}(X_t^x) D$.

\noindent Hence, arguing as in Proposition \ref{main}, we deduce that $\displaystyle |Z^{\alpha,x,n}_t| \leq |\nabla _x v^{\alpha, n} (X^{x}_t) G(X^{x}_t)| \leq {C}/{\mu}$ and $ \displaystyle |U^{\alpha,x,n}_t| \leq |\nabla _x v^{\alpha, n} (X^{x}_t)D| \leq \frac{C}{\mu}$, with $C$ independent of $n$ and $\alpha$.

\noindent Moreover, see \cite[Prop 5.2]{FuTessAOP2002},  the map 
$ x \to (Y^{x,\alpha,n}_t, Z^{x,\alpha,n}_t,  U^{x,\alpha,n}_t)$  is Gateaux differentiable and  the equation for the derivative in the direction $h \in H$, $|h|=1$,  is the following:
\begin{align*}
\nabla _xY^{x,\alpha,n}_t h   &= \nonumber
  \int_t^n [\phi ^{h, \alpha}(s, \nabla _xZ ^{x,\alpha,n}_sh, \nabla _xU^{x,\alpha,n}_sh) - \alpha \nabla_ x Y^{x,\alpha,n}_sh]  \, ds - \int_t^n \nabla_x  Z^{x,\alpha,n}_s h \, dW^1_s \\ & - \int_t^n \nabla_x  U^{x,\alpha,n}_s h \, dW^2_s, \qquad 0 \leq t \leq n.
\end{align*}
where 
\begin{align*} &  \phi^{h,\alpha,n}  (s,z,u)=   \nabla _x \psi(X^{x}_s, Z^{x,\alpha,n}_sG^{-1}(X^{x}_s),  U^{x,\alpha,n}_s)
\nabla_x X^x_s h   +  \nabla _u \psi(X^{x}_s, Z^{x,\alpha,n}_sG^{-1}(X^{x}_s), U^{x,\alpha,n}_s) u h  \\& +  \nabla _z \psi(X^{x}_s, Z^{x,\alpha,n}_sG^{-1}(X^{x}_s), U^{x,\alpha,n}_s)[
Z^{x,\alpha,n}_s   \nabla_ x G^{-1}(X^x_s) \nabla_x X^x_s h +  
z \, h  G^{-1}(X^x_s) ]\end{align*}
Notice that $\phi^{h,\alpha} (t,z,u)$ is affine in $z$ and $u$ and :
\[   | \phi^{h,\alpha,n} (s,z,u)- \phi^{h,\alpha,n} (s,0,0)| \leq L_u |u|+ L_z  M_{G^{-1}} |z|\leq (L_z^2 M^2_{G^{-1}}+ L^2_u)^{1/2} ( |z|^2 + |u|^2)^{1/2}, \qquad \mathbb{P}-a.s.\]
where here  and in the following the constant $C$ may change from line to line  but always  independently from $n$, $\varepsilon$ and from $\alpha$.

We can apply Lemma \ref{ParRasFin} with $\nu_s =  (L_z^2 M^2_{G^{-1}}+ L^2_u)  s = : K s$, indeed for  $  \e = \frac{1}{2}(\mu -  2K)$, we have, recalling also that $ U^{x,\alpha,n}_s $ and $ Z^{x,\alpha,n}_s$ are bounded uniformly in $s$, $\alpha$ and $n$
\begin{equation}\label{stima}
\E \left[\int^n_0 | \phi^{h,\alpha,n}(s,0,0)| e^{(-\alpha+K) s  }  \, dt  \right]^2 \leq \frac{C}{\varepsilon} \int^n_0  e^{(\varepsilon-2\alpha +2K) s}\E |  \nabla_x X^x_s h | ^2   \, dt   \leq   \frac{ C } { \mu- 2 K }.
\end{equation}
Therefore the following estimate holds, arguing as before in \eqref{stima}, for all $0\leq t\leq n$: 
\begin{multline}\label{stimahofinito}
\E  \sup_{s \in [t,n]}e^{2(-\alpha+K) s } |\nabla _x Y^{x,\alpha,n}_s h|^2    + \E \int_t^n \! e^{2(-\alpha+K) s } \left[  |\nabla _x Z^{x,\alpha,n}_s h |^2 + |\nabla _x U^{x,\alpha,n}_s h |^2\right] \, dt  \\ \leq C\,  \E \left[  \int_t^n e^{ (-\alpha+K) s}   |\phi^{h,\alpha,n}(s,0,0)| \, ds \right] ^2  \leq  \frac{ C e^{( -2\alpha - \frac{1}{2}\mu +K  ) t }} { \mu- 2 K} , \quad    t \leq s \leq n.
\end{multline}
In particular, we have for all $t\geq 0$:
\begin{equation}\label{stima0}
\E \Big( \,  e^{  2 K t} |\nabla _x Y^{x,\alpha,n}_t h|^2  \Big) \leq   {C }\,  e^{(-\frac{1}{2}\mu +  K )  \, t  }.
\end{equation}
From estimate  \eqref{stimahofinito}  we deduce that $(\nabla _x Y^{x,\alpha,n} h,\nabla _x Z^{x,\alpha,n} h, \nabla _x U^{x,\alpha,n} h)$ weakly converges in  the Hilbert space  $ L^2 (\Omega \times (0,T); \mathbb{R} \times \Xi^* \times H^*)$  to some $(R^{x,\alpha,h}, V^{x,\alpha,h},  M^{x,\alpha,h})$, for every $T>0$. 
 From \eqref{stima0} we also  
 have that $\nabla _x Y^{x,\alpha,n}_0 h$  converge in $\R$ to $\xi^{x,\alpha,h}$. 

We define for every $t\geq 0$
$$
\tilde{R}^{x,\alpha,h}_t  =  \xi^{x,\alpha,h}+ 
  \int_0^t \left[\phi^{h,\alpha}(s, V^{x,\alpha,h}_s,M^{x,\alpha,h}_s) - \alpha R^{x,\alpha,h}_s \right] ds \,- \int_0^t V^{x,\alpha,h}_s \, dW^1_s 
- \int_0^t M^{x,\alpha,h}_s \, dW^2_s.
$$
Now we compare the above with the forward equation fulfilled by
 $(\nabla _x Y^{x,\alpha,n} h,\nabla _x Z^{x,\alpha,n} h, \nabla _x U^{x,\alpha,n} h)$, namely:
\begin{align*}
\nabla _xY^{x,\alpha,n}_t h  =&  \nabla _xY^{x,\alpha,n}_0 h
 +  \int_0^t \left[\phi^{h,\alpha,n}(s,  \nabla_x  Z^{x,\alpha,n}_s, \nabla_x  U^{x,\alpha,n}_s) - \alpha  \nabla_x  Y^{x,\alpha,n}_s  h \right] ds \\ &- \int_0^t  \nabla_x  Z^{x,\alpha,n}_s h\, dW^1_s 
- \int_0^t  \nabla_x  U ^{x,\alpha,n}_s h \, dW^2_s, \qquad  \mathbb{P}-a.s..
\end{align*}

Since every term in the R.H.S., passing to a subsequence if necessary, weakly converges in $L^2 (\Omega \times (0,T); \mathbb{R})$,  see also \cite[Theo. 3.1]{HuTess}, we have that $\tilde{R}^{x,\alpha,h}_t= {R}^{x,\alpha,h}_t,$ \,$ \mathbb{P}-$a.s. for a.e. $t \geq 0$.
 Thus the triplet processes $(\tilde{R}^{x,\alpha,h}, V^{x,\alpha,h},  M^{x,\alpha,h})$ verifies for all $t>0$, $\mathbb{P}$-a.s.:
$$
\tilde{R}^{x,\alpha,h}_t  =  \tilde{R}^{x,\alpha,h}_0 + 
  \int_0^t \left[\phi^{h,\alpha}(s, V^{x,\alpha,h}_s,M^{x,\alpha,h}_s) - \alpha \tilde{R}^{x,\alpha,h}_s \right] ds \,- \int_0^t V^{x,\alpha,h}_s \, dW^1_s 
- \int_0^t M^{x,\alpha,h}_s \, dW^2_s.
$$
 where 
\begin{align*} &  \phi^{h,\alpha}  (s,z,u)=   \nabla _x \psi(X^{x}_s, Z^{x,\alpha}_sG^{-1}(X^{x}_s),  U^{x,\alpha}_s)
\nabla_x X^x_s h   +  \nabla _u \psi(X^{x}_s, Z^{x,\alpha}_sG^{-1}(X^{x}_s), U^{x,\alpha}_s) u h  \\& +  \nabla _z \psi(X^{x}_s, Z^{x,\alpha}_sG^{-1}(X^{x}_s), U^{x,\alpha}_s)[
Z^{x,\alpha}_s   \nabla_ x G^{-1}(X^x_s) \nabla_x X^x_s h +  
z \, h  G^{-1}(X^x_s) ]\end{align*}

Moreover, thanks to  \eqref{stimahofinito} and \eqref{stima0} we have that 
\begin{equation}\label{stimahofinito2tilde}
 \E \sup_{s \in [0,T]}e^{ 2 K s }|\tilde{R}^{x,\alpha,h}_s|^2  < + \infty \quad\text{ and }\quad \E \, e^{2  K s} |\tilde{R}^{x,\alpha,h}_s |^2  \leq   \tilde{C }\,  e^{(-\mu + 2 K )  \, s  }, 
\end{equation}
therefore, $(\tilde{R}^{x,\alpha,h}, V^{x,\alpha,h}, M^{x,\alpha,h})$ is
the unique solution of equation:
\begin{equation}\label{eq-inf-hor-R}d_sR_s=[\phi^{h,\alpha}(s,V_s,M_s)-\alpha R_s]ds -V_sdW^1_s-M_sdW^2_s\end{equation} 
in the class of processes with the regularity imposed in Lemma \ref{ParRasInf} veryfying:
\begin{equation}\label{cond}
 \E \sup_{t \in [0,T]}|\tilde{R}^{x,\alpha,h}_t|^2  < + \infty \quad\text{ and }\quad  \lim_{T \to +\infty} \E\,  e^{ 2K2 T}|\tilde{R}^{x,\alpha,h}_T|^2 =0, \qquad \forall T >0.
\end{equation}
We then closely follow the proof of \cite[Prop 3.2]{HuTess}, indeed we get that
$\lim_{n \to + \infty}\nabla_x Y^{\alpha, n,x}_0 h=\tilde{R}^{\alpha, x,h}(0)$,  defines a linear and bounded operator $\tilde{R}^{\alpha,x}(0)$ from $H$ to $H$, by \eqref{stima0},   such that  $\tilde{R}^{\alpha,x}(0)h=\tilde{R}^{x,\alpha,h}(0)$, moreover for every fixed $h \in H$, $x \to \tilde{R}^{\alpha,x}(0)h$ is continuous in $x$, we will sketch the argument by the the end of the proof in a similar point. Therefore, by dominated convergence, we get that:
\begin{multline}\label{difffin}
\lim_{\ell \downarrow 0} \frac{Y^{x+\ell h,\alpha}_0- Y^{x,\alpha}_0}{\ell }=\lim_{\ell \downarrow 0} \lim_{n\to \infty}\frac{Y^{x+\ell h,\alpha, n}_0- Y^{x,\alpha, n}_0}{\ell}= \lim_{\ell \downarrow 0} \lim_{n\to \infty} \int_0^1 \nabla_ x Y^{,x+\theta \ell h,\alpha, n}_0 h \, d\theta\\=  \lim_{\ell  \downarrow 0}  \int_0^1  \tilde{R}^{ x+\theta \ell h,\alpha}(0)h\, d\theta =  \tilde{R}^{x,\alpha}(0)h.
\end{multline}
Thus $v^{\alpha }																		$ is differentiable and since $Y^{x,\alpha}_t=v^{\alpha }(X^x_t)	$ we have
 $\nabla _xY^{x,\alpha}_t h =v^{\alpha}(X^x_t)\nabla_x X^x_th$. 
 
Fixing $T>0$ we can see the equation satisfied by 
$(Y^{x,\alpha},Z^{x,\alpha},U^{x,\alpha})$ as a BSDE on  $[0,T]$ with final condition $v^{\alpha}(X^x_T)$  and  we can apply standard results on the differentiability of  markovian, finite horizon BSDEs (see, for instance, \cite{FuTessAOP2002}) to deduce that the map $ x\to Y^{x,\alpha}$ is of class $\mathcal{G}^1$ from $H$ to $ L^2 _ {\mathcal{P}}(\Omega,; C([0,T]; \mathbb{R}))$ and $ x \to Z^{x,\alpha} $ is of class $\mathcal{G}^1$ from $L^2_{\mathcal {P}}([0,T] \times \Omega; \Xi^*)$.
Moreover for every $h \in  H$, for every $0\leq t\leq T$ it holds that:
\begin{align}\label{nablaYalpha}
\nabla _xY^{x,\alpha}_t h   &= \nabla_x  Y^{x,\alpha}_T h +\nonumber
  \int_t^T [\phi ^h(s, \nabla _xZ ^{x,\alpha}_sh, \nabla _xU^{x,\alpha}_sh) - \alpha \nabla_ x Y^{x,\alpha}_sh]  \, ds \\ & - \int_t^T \nabla_x  Z^{x,\alpha}_s h \, dW^1_s  - \int_t^T \nabla_x  U^{x,\alpha}_s h \, dW^2_s, \qquad 0 \leq t \leq n.
\end{align}
Comparing the above with \eqref{eq-inf-hor-R} and noticing that for all $T>0$: 
$$\mathbb{E} e^{2 K T}|\nabla _xY^{x,\alpha}_T h|^2
=\mathbb{E} e^{2 K T}| \nabla_x v^{\alpha}(X^x_T)\nabla_x X^x_T h|^2\leq C e^{(2 K-\mu) T} 	$$
the uniqueness part of Lemma \ref{ParRasInf} tells us that $(\nabla _xY^{x,\alpha}_{\cdot} h,
\nabla _xZ^{x,\alpha}_{\cdot} h, \nabla _xU^{x,\alpha}_{\cdot} h )$ coincides with $(\tilde{R}^{x,h,\alpha},
V^{x,h,\alpha},M^{x,h,\alpha})$ and is the unique solution of equation  \eqref{eq-inf-hor-R}  in the sense of Lemma \ref{ParRasInf}.

\medskip 
\noindent\textit{Part II - Differentiability of $\bar{v}$}

\noindent

We also introduce the following infinite horizon BSDE:
\begin{equation}\label{Ualphadiffrho0}
- d\,{R}^{x,h}_s  = \phi^h(s,V^{x,h}_s,M^{x,h}_s)ds- V^{x,h}_t \, dW_t^1  - M^{x,h}_t \, dW_t^2 \qquad t \geq 0.
\end{equation}
with 
\begin{align*} \phi^h(s,z,u)= &
 [ \nabla _x \psi(X^{x}_s, \bar{Z}^{x}_sG^{-1}(X^{x}_s), \bar{U}^{x}_s)+ 
 \nabla _z \psi(X^{x}_s, \bar{Z}^{x}_sG^{-1}(X^{x}_s), \bar{U}^{x}_s)
\bar{Z}^{x}_s   \nabla_ x G^{-1}(X^x_s)] \nabla_x  X^x_s h
   \\&  + \nabla _u \psi(X^{x}_s, \bar{Z}^{x}_sG^{-1}(X^{x}_s), \bar{U}^{x}_s)
u  + \nabla _z \psi(X^{x}_s, \bar{Z}^{x}_sG^{-1}(X^{x}_s), \bar{U}^{x}_s)
z
 \end{align*}
By Lemma \ref{ParRasInf} has a unique solution in the class  of processes ${R}^{x,h} \in L^{2, loc}_\Pc (\Omega;{C} ([0,+\infty[;\mathbb{R}))$, ${V}^{x,h}\in L^{2,loc}_\Pc (\Omega \times [0,+\infty [;\Xi ^*) $, $M\in
L^{2,loc}_\Pc (\Omega \times [0,+\infty [;H ^*)$ verifying:
\begin{equation}\label{cond1}
\lim_{T\to +\infty} e^{ 2 K T} \E\,  |{R}^{x,h}_T|^2 =0, \qquad \forall T >0.
\end{equation}
As in \cite[Theorem 5.1]{FuhHuTess} we claim that, along the sequence  $( \alpha_m)$ introduced in \eqref{defvbar},  it holds:
\begin{equation}\label{convdiff}
\nabla _x {v}^{\alpha_m}(x) h=\nabla _x Y^{\alpha_m	,x}_0h = R^{x,\alpha_m,h}_0 \to R^{x,h}_0,
\end{equation}
as $m\rightarrow \infty$.

Let us introduce again some  parabolic approximations: for:
\begin{equation*} \left \{
\begin{array}{ll}
- d\,{R}^{x,\alpha, n,h}_s\! = &\!\! \phi^{h,\alpha}(s,V^{x,\alpha, n,h}_s,M^{x,\alpha, n,h})ds-\alpha {R}^{x,\alpha, n,h}_s \, ds - V^{x,\alpha, n,h}_s \, dW_s^1  - M^{x,\alpha, n,h}_s \, dW_s^2 \quad s \in [0,n], \\
\ \ R^{x,\alpha, n,h}_n\! =& \!\!0
\end{array} \right.
\end{equation*}
and 
\begin{equation*} \left \{
\begin{array}{ll}
- d\,{R}^{x,n,h}_s & =  \phi^{h}(s,V^{ x,h,n}_s,M^{x,n,h})ds- V^{ x,h,n}_s \, dW_s^1  - M^{x,h}_s \, dW_s^2 \qquad s \in [0,n], \\
\ \ \ R^{ x,h,n}_n&=0
\end{array} \right.
\end{equation*}
Since along the sequence $(\alpha_m)$ selected in Section \ref{backward} we have
$$\mathbb{E}\sup_{s\in [0,n]} |\bar{Y}^{x}_s-{Y}^{x,\alpha_m}_s|^2+\mathbb{E}\int_0^n\left[ |\bar{Z}_s-{Z}^{x,\alpha_m }_s|^2+ |\bar{U}^{x}_s-U^{x,\alpha_m }_s|^2 \right] ds\rightarrow 0$$
and consequently
$$\mathbb{E}\int_0^n |\phi^{h,\alpha_m}(s,0,0)-\phi^{h}(s,0,0)|^2 ds\rightarrow 0\quad \hbox{as $m\rightarrow \infty$}$$
standard estimates on finite horizon BSDEs give:
\begin{equation}\label{limiteUN}
\mathbb{E}\sup_{s\in [0,n]} |{R}^{x,n,h}_s-{R}^{x,\alpha_m, n,h}_s|^2\rightarrow 0,\quad \hbox{as $m\rightarrow \infty$}.\end{equation}
 Moreover if we compare with the solution $(\tilde{R}^{x,\alpha,h}, V^{x,\alpha,h}, M^{x,\alpha,h})$  of equation \eqref{eq-inf-hor-R}
\begin{equation} \left \{
\begin{array}{ll}
-d\,({R}^{x,\alpha,n,h}_s\!\! -\tilde{R}^{x,\alpha,h}_s )& \!\!\!\!= \phi^{h,\alpha}(s,{V}^{x,\alpha,n,h}_s \!-{V}^{x,\alpha,h}_s, {M}^{x,\alpha,n,h}_s \!-M^{x,\alpha,h}_s )ds
   -\alpha [{R}^{x,\alpha,n,h}_s -\tilde{R}^{x,\alpha,h}_s ] \, ds\\&-[{V}^{x,\alpha,n,h}_s -{V}^{x,\alpha,h}_s] \, dW_s^1-[{M}^{x,\alpha,n,h}_s -{M}^{x,\alpha,h}_s] \, dW_s^2,  \\
\ \ \ \ \  {R}^{x,\alpha,n,h}_n -\tilde{R}^{x,\alpha,h}_n &\!\!\!\!= - \nabla _x  {v}^\alpha (X^x_n) \nabla_x  X^x_n h
\end{array} \right.
\end{equation}
Thus Lemma  \ref{ParRasFin} estimate \eqref{stimaparrasfin} yields:
\begin{equation}\label{asintalpha}
| R^{x,\alpha,n,h}_0- \tilde{R}^{x,\alpha,h}_0 |^2\leq 
\E \left( e^{2 kn} |\nabla _x  {v}^\alpha (X^x_n) \nabla_x  X^x_n h |^2 \right) \leq  C  e^{(2K- \mu) n} \to 0,  \ \text{ as } n \to +\infty.
\end{equation}
Notice that the right hand side does not depend on $\alpha$ . 
Finally
\begin{equation} \left \{
\begin{array}{ll}
-d\,({R}^{x,n,h}_s -{R}^{x,h}_s )& = \phi^{h}(s,{V}^{x,n,h}_s -{V}^{x,h}_s, {M}^{x,n,h}_s -{M}^{x,h}_s )ds
    \\&-[{V}^{x,n,h}_s -{V}^{x,h}_s] \, dW_s^1-[{M}^{x,n,h}_s -{M}^{x, h}_s] \, dW_s^2,  \\
\ \ \ \ \  {R}^{x,n,h}_n -{R}^{x,h}_n &= - \tilde{R}^{x,h}_n
\end{array} \right.
\end{equation}
and taking into account  \eqref{cond1}, one has, again by Lemma \ref{ParRasFin} relation \eqref{stimaparrasfin}:
\begin{equation}\label{asint0}
 | R^{x,n,h}_0- {R}^{x,h}_0 |^2\leq 
\E \left( e^{2 K n } | R^{x,h}_n| ^2 \right)\leq C  
 e^{(2 K -\mu)n } \to 0, \qquad \text{ as } N \to +\infty.
\end{equation}

Therefore summing up \eqref{asintalpha}, \eqref{asint0} and \eqref{limiteUN} we have that:
\[ R^{x,\alpha_m,h}_0 \to R^{x,h}_0, \qquad \text{ as } m \to +\infty.\]
Finally the continuity with respect to $x$ of $R^{x,h}_0$ descends immediately from  \eqref{asint0} and from the continuity of the map  $ x \to R_0^{x,n,h}$ proved in \cite[Prop. 4.3]{FuTessAOP2002}.

\medskip
We can now conclude as above (and ass in \cite[Prop 3.2]{HuTess});
${R}^{x,h}(0)$,  defines a linear and bounded operator ${R}^{x}(0)$ from $H$ to $H$, such that  ${R}^{x}(0)h={R}^{x,h}(0)$, and  we  have:

\begin{align*}
\lim_{t \downarrow 0} \frac{\bar{v}(x+th)-\bar{v}(x)}{t}
=\lim_{t \downarrow 0}\frac{\bar{Y}^{x+th}_0- \bar{Y}^{x}_0}{t} =
\lim_{t \downarrow 0} \lim_{m \to 0}\frac{{Y}^{x+th,\alpha_m}_0- {Y}^{x,\alpha}_0}{t} =\\= \lim_{t \downarrow 0} \lim_{m \to 0} \int_0^1 \nabla_ x Y^{x+\theta th,\alpha_m}_0 h \, d\theta= \lim_{t \downarrow 0} \lim_{m \to 0} \int_0^1   R^{ x+\theta th,\alpha_m,h}(0)h\, d\theta = \\= \lim_{t \downarrow 0}  \int_0^1  R^{ x+\theta th}(0) h\, d\theta =  {R}^{x}(0)h.
\end{align*}

\ep

\bibliographystyle{plain}
\bibliography{biblionew}

\end{document}